\documentclass[12pt]{article}
\usepackage{amssymb,amsfonts,amsmath, psfrag,eepic,colordvi,graphicx,epsfig,ytableau}
\usepackage{amssymb,latexsym,graphics,array}
\usepackage{MnSymbol}
\usepackage[enableskew]{youngtab}
\usepackage{float,enumerate,enumitem}
\usepackage{tikz,ifthen}
\usepackage{pgfplots}
\usetikzlibrary{math}
\topskip  
\newcommand{\rmnum}[1]{\romannumeral #1}

\numberwithin{equation}{section}
\newtheorem{theorem}{Theorem}[section]

\newtheorem{remark}[theorem]{Remark}
\newtheorem{lemma}[theorem]{Lemma}

\usepackage{makecell} 
\usepackage[numbers,sort&compress]{natbib}

\begin{document}
	\parskip 6pt
 
	\bibliography{ref}
	\bibliographystyle{refstyle}
	\pagenumbering{arabic}
	\def\sof{\hfill\rule{2mm}{2mm}}
	\def\ls{\leq}
	\def\gs{\geq}
	\def\SS{\mathcal S}
	\def\qq{{\bold q}}
	\def\MM{\mathcal M}
	\def\TT{\mathcal T}
	\def\EE{\mathcal E}
	\def\lsp{\mbox{lsp}}
	\def\rsp{\mbox{rsp}}
	\def\pf{\noindent {\it Proof.} }
	\def\mp{\mbox{pyramid}}
	\def\mb{\mbox{block}}
	\def\mc{\mbox{cross}}
	\def\qed{\hfill \rule{4pt}{7pt}}
	\def\block{\hfill \rule{5pt}{5pt}}
	\def\lr#1{\multicolumn{1}{|@{\hspace{.6ex}}c@{\hspace{.6ex}}|}{\raisebox{-.3ex}{$#1$}}}
	\renewcommand{\figurename}{Fig.}
	\def\blue{\textcolor{blue}}
	\def\red{\textcolor{red}}
	\def\green{\textcolor{green}}
	\def\pink{\textcolor{pink}}
	\def\violet{\textcolor{violet}}
	\def\cyan{\textcolor{cyan}}
	\def\magenta{\textcolor{magenta}}
	\def\purple{\textcolor{purple}}

	\begin{center}
		{\Large \bf   A Self-Conjugate Partition Analog of $(t,t+1)$-Core Partitions with Distinct Parts}
	\end{center}
	
	\begin{center}
 
{\small  Huan Xiong  
			and Lihong Yang$^{*}$\footnote{$^*$Corresponding author.}  }
\vskip 2mm

 Institute for Advanced Study in Mathematics \\
   Harbin Institute of Technology\\
   Heilongjiang 150001, P.R. China\\[6pt]

   \vskip 2mm

Email:     
huan.xiong.math@gmail.com, 
yanglihong954@163.com 

\end{center}
	
	\noindent {\bf Abstract.}  
	Simultaneous core partitions have been widely studied in the past 20 years. In 2013, Amdeberhan gave several conjectures on the number, the average size, and the largest size of $(t,t+1)$-core partitions with distinct parts, which was proved and generalized by Straub, Xiong, Nath-Sellers, Zaleski-Zeilberger, Paramonov, and many other mathematicians. In this paper, we introduce a proper self-conjugate partition analog of $(t,t+1)$-core partitions with distinct parts, and derive the number, the average size, and the largest size for such core partitions.

% The generating function of $t$-core partitions with distinct parts is obtained. We also prove the results on the number, the largest size, and the average size of $(t, t + 1)$-core partitions. This gives a complete answer to a conjecture of Amdeberhan, which is partly and independently proved by Straub, Nath and Sellers, and Zaleski recently.


	
	%===========================================================================
	
	\section{Introduction}
	
A partition $\lambda$ of a positive integer $n$ is a finite non-increasing sequence of positive integers $( \lambda_1,\lambda_2,\dots,\lambda_m )$ satisfying $\lambda_1+\lambda_2+\dots+\lambda_m=n$ (see \cite{Macdonald, ec2}). We denote it by $\lambda \vdash n$. We say that  $\lambda_i ~ (1\leq i \leq m)$ are the parts of $\lambda$ and $|\lambda|:=\sum_{1\leq i\leq m}\lambda_i$ is the size of $\lambda$.  The number $m$ of parts of $\lambda$ is called the length of $\lambda$, which is denoted by $\ell(\lambda)$.
The {\em Young diagram} of $\lambda$ is defined to be an up- and left-justified array of $n$ boxes with $\lambda_i$ boxes in the $i$-th row. The {\em conjugate partition} of $\lambda$ is the partition whose Young diagram is obtained by reflecting the Young diagram of $\lambda$ about the main diagonal, and $\lambda$ is said to be  {\em self-conjugate} it is equal to its conjugate partition. For the $(i,j)$-box in the $i$-th row and $j$-th column in the  Young diagram, its {\em hook length} is the number of boxes directly to the right, and directly to the below, including the boxes itself. For example, the Young diagram and the hook lengths of a self-conjugate partition $(6,4,2,2,1,1)$ are given in Figure \ref{fig:YD}.
	
		\begin{figure} [h]
		\centering
		{$\begin{array}[b]{*{6}c}\cline{1-6}
				\lr{11}&\lr{8}&\lr{5}&\lr{4}&\lr{2}&\lr{1}\\\cline{1-6}
				\lr{8}&\lr{5}&\lr{2}&\lr{1}\\\cline{1-4}
				\lr{5}&\lr{2}\\\cline{1-2}
				\lr{4}&\lr{1}\\\cline{1-2}
				\lr{2}\\\cline{1-1}
				\lr{1}\\\cline{1-1}
			\end{array}$}
		\caption{The Young diagram and the hook lengths of a self-conjugate partition $(6,4,2,2,1,1)$.}\label{fig:YD}
	\end{figure}
	
	For a positive integer $t$, a partition is considered a {\em $t$-core partition} if none of its hook lengths is divisible by $t$.  Let $t_i~(1\leq i\leq k)$ be some positive integers. We say that $\lambda$ is a \emph{$(t_1,t_2,\ldots, t_k)$}-core partition if it is simultaneously a $t_1$-core, a $t_2$-core, $\ldots$, and a $t_k$-core partition (see \cite{tamd, KN}). For instance, we can see from Figure \ref{fig:YD} that $(6,4,2,2,1,1)$ is a $(6,7,13)$-core partition.    
 Several statistics of simultaneous core partitions, such as the numbers, the largest sizes and the average sizes, have been widely studied in the past twenty years (see  \cite{tamd1,  and, AHJ, CHW, HJJ,HJJ1,HK,CEZ, ford,HHL, PJ,NK,WX,HM,N3,NS,  SNKC,   SZ,Wang,Xiong2, Xiong3, Xiong4,  SDH, SYH,  Za,Za2,ZZ}). For example, when $s$ and $t$ are coprime, Anderson \cite{and} showed that the number of $(s,t)$-core partitions equals the rational Catalan number $\frac{1}{s+t}{s+t \choose s}$. In 2014, Armstrong, Hanusa, and Jones \cite{AHJ} conjectured that the average size of $(s,t)$-core partitions equals$\frac{(s+t-1)(s-1)(t-1)}{24}$ when $s$ and $t$ are coprime, which was proved by Johnson \cite{PJ} and Wang \cite{Wang} independently.

A partition is a self-conjugate partition if its Young diagram is symmetric along the main diagonal. When $s$ and $t$ are coprime, Ford, Mai, and Sze \cite{ford} derived the number of self-conjugate $(s, t)$-core partition, which is equal to $\binom{\lfloor\frac{s}{2}\rfloor+\lfloor\frac{t}{2}\rfloor}{\lfloor\frac{s}{2}\rfloor}$; while the average size of such partition was conjectured to be $\frac{(s+t-1)(s-1)(t-1)}{24}$ by Armstrong, Hanusa, and Jones \cite{AHJ} in 2014 and proved by Chen, Huang, and Wang \cite{CHW} in 2016. Furthermore, Olsson and Stanton \cite{ols} proved that the largest size for such partition is equal to $\frac{(s^2-1)(t^2-1)}{24}$.

In 2013, Amdeberhan \cite{tamd} provided the following conjectures on the number, the average size, and the largest size of $(t,t+1)$-core partitions with distinct parts, which was proved and generalized by Straub \cite{Straub}, Xiong \cite{Xiong2}, Nath-Sellers \cite{NS}, Zaleski-Zeilberger \cite{Za, ZZ}, Paramonov \cite{Para} and many other mathematicians.

\begin{theorem}{\upshape   (see \cite{tamd,Straub,Xiong2})}\label{conj:amd_conjecture}
Let $t\geq 1$ be a given positive integer and $(F_i)_{i\geq 1}=(1,1,2,3,5,8,13,\ldots)$ be the  Fibonacci numbers.  We have the following results for $(t, t + 1)$-core partitions with distinct parts.
\begin{enumerate}
\item[\upshape (1)] The number of such partitions is  $F_{t+1}$, which means that $F_t=F_{t-1}+F_{t-2}$ for $t\geq 3$ and it has the ordinary generating function
\begin{equation}\label{egen1*}
	\sum_{t\geq0}F_{t+1}x^t=\frac{1}{1-x-x^2}.
\end{equation}

\item[\upshape (2)]The total sum of the sizes of these partitions and the average size are, respectively, given by
$$\sum\limits_{\substack{i+j+k=t+1\\ i,j,k\geq 1}}F_iF_jF_k\ \ \ \text{and}\ \  \ \sum\limits_{\substack{i+j+k=t+1\\ i,j,k\geq 1}}\frac{F_iF_jF_k}{F_{t+1}}.$$
%which has the ordinary generating function
%\begin{equation}\label{egen3*}
%\frac{-6+x+11x^2+6x^3}{50(1-x-x^3)}.
%\end{equation}

\item[\upshape (3)] The largest size of such partitions is
$\lfloor\frac13\binom{t+1}{2}\rfloor$, where  $\lfloor x \rfloor$ denotes the largest integer not greater than $x$.

\item[\upshape (4)] The number of such partitions with the largest size is $2$ if $t  \equiv 1 \, (\text{mod} \ 3)$ and $1$ otherwise.
\end{enumerate}
\end{theorem}

This paper aims to derive a proper self-conjugate partition analog of Theorem~\ref{conj:amd_conjecture}. To achieve this, one may start by considering self-conjugate $(t,t+1)$-core partitions with distinct parts. However, the set of such partitions has a straightforward structure, which is equal to $\{ (n,n-1,n-2,\ldots, 3,2,1): n\leq t/2 
 \}$ (see \cite{NS}), thus the number, the average size, and the largest size are easy to determine. Therefore, such partitions may not be a proper analog of $(t,t+1)$-core partitions with distinct parts.  Instead, in this paper, we introduce and study the set $\mathcal{DS}(t)$ of self-conjugate $(t,t+1)$-core partitions whose first $s(\lambda)$ parts are distinct, where $s(\lambda)$ is the size of the Durfee square \cite{Andrews} of $\lambda$.
 Recall that a partition $\lambda$ has a Durfee square of size~$s$ if $s$ is the largest integer such that the partition contains at least $s$ parts which are larger than or equal to~$s$ (see \cite{Andrews}).
We always assume that the empty partition $\emptyset$ belongs to $\mathcal{DS}(t)$. 
 Actually, studying such partitions is reasonable since a self-conjugate partition is uniquely determined by its first $s(\lambda)$ parts. For example, the partition $(6,4,2,2,1,1) \in \mathcal{DS}(6)$ since it is a $(6,7)$-core, $s(\lambda)=2$, and the first two parts $6$ and $4$ are distinct.  
 In Table \ref{table:set}, we list the number of partitions in $\mathcal{DS}(t)$ and the set of such partitions for $1 \leq t \leq 7$. 
 We derive the following results for the partitions in~$\mathcal{DS}(t)$.

\begin{table}[H]
	\centering
	\scriptsize
	\renewcommand\arraystretch{1.2}
	\caption{The number of partitions in $\mathcal{DS}(t)$ and the set of such partitions.}\label{table:set}
	\vskip 2mm
	\begin{tabular}{|m{0.8cm}<{\centering}|c|c|	}	\hline 
		$t$ &   $The \,\, set\, ~\mathcal{DS}(t)$    &  $ The \,\, number \,f(t) $ \\ \hline
		
		\vskip 1.5mm
       $1$  & $\{\emptyset\}$ & $1$\\ \hline	
       \vskip 1.5mm
      $2$  & $\{\emptyset, (1)\}$ & $2$\\ \hline	

       \vskip 1.5mm
        $3$  & $\{\emptyset, (1), (3,1,1)\}$ & $3$\\ \hline	

        \vskip 1.5mm
        $4$  & $\{\emptyset, (1), (2,1), (4,1,1,1)\}$ & $4$\\ \hline	

        \vskip 1.5mm
        $5$  & $\{\emptyset, (1), (2,1), (4,1,1,1), (5,1,1,1,1), (4,2,1,1)\}$ & $6$\\ \hline

        \vskip 1.5mm
        $6$  & $\thead{ \{\emptyset, (1), (2,1), (3,1,1),  (3,2,1), (5,1,1,1,1), (5,2,1,1,1), \\  (6,1,1,1,1,1), (6,4,2,2,1,1)\}}$ & $9$\\ \hline

        \vskip 1.5mm
        $7$  & $\thead{\{\emptyset, (1), (2,1), (3,1,1), (3,2,1), (5,1,1,1,1),  (5,2,1,1,1), \\ (5,3,2,1,1), (6,1,1,1,1,1), (6,2,1,1,1,1), (7,1,1,1,1,1,1), \\  (7,4,2,2,1,1,1), (7,6,2,2,2,2,1)\}}$ & $13$\\ \hline

	\end{tabular}
  \end{table}

 \begin{theorem}\label{gen1} 
		Let $t \geq 1$ be a positive integer. We have the following results for self-conjugate partitions in $\mathcal{DS}(t)$.
		\begin{enumerate}
\item[\upshape (1)]  Let $f(t)$ be the number of partitions in $\mathcal{DS}(t)$ and set $f(-1)=f(0)=1$. Then we have $f(1)=1$,   $f(2)=2$, and for $t\geq 3$, we have $f(t)=f(t-1)+f(t-3)$. Thus the ordinary generating function of $f(t-1)$ is 
\begin{equation}\label{egen1}
	\sum_{t\geq0}f(t-1)x^t=\frac{1}{1-x-x^3}.
\end{equation}

\item[\upshape (2)]   Let $g(t)$ be the sum of the number of main diagonal hooks of all partitions in $\mathcal{DS}(t)$ and set $g(-1)=g(0)=0$. Then $g(1)=0$, $g(2)=1$, and the   ordinary generating function of $g(t-1)$ is 
{\begin{equation}\label{egen2}
	\sum_{t\geq0}g(t-1)x^t=\frac{x^3}{(1-x-x^3)^2}.
\end{equation}}

\item[\upshape (3)]   Let $h(t)$ be the sum of the sizes of all partitions in $\mathcal{DS}(t)$ and  set $h(-1)=h(0)=0$. Then $h(1)=0$, $h(2)=1$, and the ordinary generating function of $h(t-1)$ is
\begin{equation}\label{egen3}
\sum_{t\geq0}h(t-1)x^t=\frac{(1+4x-x^2+6x^3+4x^4+5x^6)x^3}{(1+x+x^3)(1-x-x^3)^3}.
\end{equation}

\end{enumerate}
	\end{theorem}

\begin{theorem}\label{largest}
Let $t$ be a positive integer and $m(t)$ be the largest size of a self-conjugate $(t,t+1)$-core partition in $\mathcal{DS}(t)$. Then we have $m(1)=0$, $m(2)=1$,   $m(t)=2t-1$ for $3 \leq t \leq 5$, and for $t \geq6$,  the following results hold.
\begin{enumerate}
\item[\upshape (1) ]  When {\upshape (\rmnum{1})} $t $ is even and $\frac{t}{2}-1-2\lfloor \frac{2t+1}{12} \rfloor \equiv 0 \, (\text{mod}\  3 )$; or {\upshape (\rmnum{2})} $t $ is odd and $\frac{t-1}{2}-2\lfloor \frac{2t+1}{12} \rfloor \equiv 2 \, (\text{mod}\  3 )$, then 
$$
m(t)=\frac{1}{3}(t-1)\left(t+1-4\lfloor \frac{2t+1}{12} \rfloor \right)+\lfloor \frac{2t+1}{12} \rfloor \left(2t+1-2\lfloor \frac{2t+1}{12} \rfloor \right);
$$
\item[\upshape (2) ]  When {\upshape (\rmnum{3})} $t $ is even and $\frac{t}{2}-1-2\lfloor \frac{2t+1}{12} \rfloor \equiv 1 \, (\text{mod}\ 3 )$; or {\upshape (\rmnum{4})} $t $ is odd and $\frac{t-1}{2}-2\lfloor \frac{2t+1}{12} \rfloor \equiv 0 \, (\text{mod}\ 3 )$, then 
$$
m(t)=\frac{1}{3}(t-1)\left(t-7-4\lfloor \frac{2t+1}{12} \rfloor \right)+\left(\lfloor \frac{2t+1}{12} \rfloor +2\right)\left(2t-3-2\lfloor \frac{2t+1}{12} \rfloor \right);
$$
\item[\upshape (3) ]  When {\upshape (\rmnum{5})} $t $ is even and $\frac{t}{2}-1-2\lfloor \frac{2t+1}{12} \rfloor \equiv 2 \, (\text{mod}\ 3 )$; or {\upshape (\rmnum{6})} $t $ is odd and $\frac{t-1}{2}-2\lfloor \frac{2t+1}{12} \rfloor \equiv 1 \, (\text{mod}\ 3 )$, then 
$$
m(t)=\frac{1}{3}(t-1)\left(t-3-4\lfloor \frac{2t+1}{12} \rfloor \right)+\left(\lfloor \frac{2t+1}{12} \rfloor+1\right) \left(2t-1-2\lfloor \frac{2t+1}{12} \rfloor \right);
$$

\item[\upshape (4)] The number of such partitions with the largest size is always equal to $1$, i.e., there always exists a unique partition in $\mathcal{DS}(t)$ with the largest size;
\end{enumerate}	
where  $\lfloor x \rfloor$ denotes the largest integer not greater than $x$.
\end{theorem}

\begin{remark}
From the above results, we can see that the number $f(t)$ of partitions in $\mathcal{DS}(t)$ satisfies the recursion $f(t)=f(t-1)+f(t-3)$, while the number $F_{t+1}$ of $(t, t + 1)$-core partitions with distinct parts satisfies the recursion $F_t=F_{t-1}+F_{t-2}$. Thus it is easy to see that the numbers of these two kinds of partitions share similar recursion formulas. Furthermore, their largest sizes and average sizes also share some similar behavior. Therefore, the partitions in $\mathcal{DS}(t)$ can be considered as a proper self-conjugate partition analog of $(t,t+1)$-core partitions with distinct parts.
\end{remark}

%In the past decade, the problem of the enumeration of simultaneous core partitions with distinct parts was raised by Amdeberhan \cite{tamd}. He conjectured explicit formulas for the number, the largest size, and the average size of $(t,t+1)$-core partitions with distinct parts,  which the first author first proved \cite{Xiong2}, and later proved independently (and extended) by Straub \cite{Straub}, Nath-Sellers \cite{NS}, Zaleski \cite{Za} and Paramonov \cite{Para}.

 For a  partition $\lambda$,  let $MD(\lambda)=\{h_1,h_2,\dots, h_k\}$ denote the set of main diagonal hook lengths, where $h_1>h_2>\cdots> h_k$ and $k=s(\lambda)$ is the number of main diagonal hooks of $\lambda$, which is also denoted by $|MD(\lambda)|$. In particular, if $\lambda=\emptyset$, then $MD(\lambda) = \emptyset$. Let $d(\lambda)$ be the sequence $(d_1,d_2,\dots,d_{k-1})$, where $d_i=h_{i}-h_{i+1}$ ($1 \leq i \leq k-1$) denotes the difference between $i$-th and $(i+1)$-th main diagonal hook lengths of $\lambda$. 
 %, that is, if $MD(\lambda)=\{h_1,h_2,\dots, h_k\}$ with $h_i>h_{i+1}$ for $i \in  [k-1]$, then $d_i=h_i-h_{i+1}$. 
 For example, let $\lambda=(5,3,2,1,1)$. Then $MD(\lambda)=\{9,3\}$ and $d(\lambda)=(d_1)=(6)$. For a nonnegative integer $m$, we say $d(\lambda)\geq m$,  if $d_i\geq m$ for all $1\leq i \leq |MD(\lambda)|-1$. Then a self-conjugate  $(t,t+1)$-core partitions  $\lambda$ belongs to the set $\mathcal{DS}(t)$ iff $d(\lambda)\geq 4$ (when $k=1$, $d(\lambda)$ is the empty sequence, thus $d(\lambda)\geq 4$ holds automatically). It is easily seen that a self-conjugate partition is uniquely determined by its main diagonal hooks, which are distinct odd numbers. In what follows, we always treat a self-conjugate partition $\lambda$ and its set of main diagonal hooks as identical. In Table \ref{table:lar}, we list the largest size of partitions in $\mathcal{DS}(t)$ and their corresponding $MD(\lambda)$ for $3\leq t \leq 26$.

\begin{table}[H]
	\centering
	\scriptsize
	\renewcommand\arraystretch{1.2}
	\caption{The largest size of partitions in $\mathcal{DS}(t)$ and their corresponding $MD(\lambda)$.}\label{table:lar}
	\vskip 2mm
	\begin{tabular}{|m{1cm}<{\centering}|c|c|m{1cm}<{\centering}|c|c|}
		\hline 
		$t$ &  $MD(\lambda)$       &  $Largest \, size$     &    $t$   &   $MD(\lambda)$     & $Largest \, size$    \\ \hline
		
		\vskip 1mm
		$3$
		&  $\{5\}$   & $5$     &  \vskip 1mm
		$15$  
		& $\{29,25,21,17\}$   & $92$                             \\ \hline
	
		\vskip 1mm
		$4$ 
		&  $\{7\}$   & $7$     &  \vskip 1mm
		$16$  
		& $\{31,27,21,15,9\}$   & $103$                             \\ \hline
		
		\vskip 1mm
		$5$ 
		&  $\{9\}$   & $9$     &  \vskip 1mm
		$17$  
		& $\{33,29,25,19,13\}$   & $119$                             \\ \hline

		\vskip 1mm
		$6$ 
		&  $\{11,5\}$   & $16$     &  \vskip 1mm
		$18$  
		& $\{35,31,27,23,17\}$   & $133$                             \\ \hline
		
		\vskip 1mm
		$7$ 
		&  $\{13,9\}$   & $22$     &  \vskip 1mm
		$19$  
		& $\{37,33,29,25,21\}$   & $145$  	                              \\ \hline
		
		\vskip 1mm
		$8$
		& $\{15,11\}$   & $26$   &  \vskip 1mm
		$20$
		&  $\{39,35,31,25,19,13\}$	   	   &  $162$                        \\ \hline
			\vskip 1mm
		$9$ 
		&  $\{17,11,5\}$   & $33$     &  \vskip 1mm
		$21$  
		& $\{41,37,33,29,23,17\}$   & $180$  	                              \\ \hline
		
		\vskip 1mm
		$10$
		& $\{19,15,9\}$   & $43$   &  \vskip 1mm
		$22$
		& $\{43,39,35,31,27,21\}$ 	  &      $196$                  \\ \hline
			\vskip 1mm
		$11$ 
		&  $\{21,17,13\}$   & $51$     &  \vskip 1mm
		$23$  
		& $\{45,41,37,31,25,19,13\}$   & $211$  	                              \\ \hline
		
		\vskip 1mm
		$12$
		& $\{23,19,15\}$   & $57$   &  \vskip 1mm
		$24$
		&  	$\{47,43,39,35,29,23,17\}$   	   &  $233$                        \\ \hline
			\vskip 1mm
		$13$ 
		&  $\{25,21,15,9\}$   & $70$     &  \vskip 1mm
		$25$  
		& $\{49,45,41,37,33,27,21\}$   & $253$  	                              \\ \hline
		
		\vskip 1mm
		$14$
		& $\{27,23,19,13\}$   & $82$   &  \vskip 1mm
		$26$
		&  	$\{51,47,43,39,35,31,25\}$  & $271$                          \\ \hline
	\end{tabular}
\end{table}

The rest of this paper is organized as follows. Section 2 proves the enumeration results of partitions in $\mathcal{DS}(t)$.
Section 3 is devoted to the proof of the largest size results for such partitions.

\section{Proof of  Theorem  \ref{gen1}}

In this section, we will give a proof of  Theorem \ref{gen1} by the recurrence relations of $f(t)$, $g(t)$, and $h(t)$.
To begin with, we need to recall the following theorem  due to
Ford, Mai and Sze \cite{ford}.

\begin{theorem}{\upshape   (see \cite{ford})}\label{pro}
Let $\lambda$ be a self-conjugate partition. Then $\lambda$ is $t$-core partition if and only if the following conditions hold.
\begin{itemize}
\item[\upshape (\rmnum{1}) ]  If $h\in MD(\lambda)$ and $h>2t$, then $h-2t\in MD(\lambda)$;
\item[\upshape (\rmnum{2}) ]  If $h_1,h_2\in MD(\lambda)$, then $h_1+h_2 \nequiv 0 \ (\text{mod}\  2t)$.
\end{itemize}
\end{theorem}

The main diagonal hook length definition lets us obtain the following result directly.

\begin{lemma}\label{sclar}
The size of a self-conjugate partition equals the sum of its main diagonal hook lengths, that is, if $\lambda$ is a self-conjugate partition, then  $$|\lambda|=\sum_{x \in MD(\lambda)}x.$$ 
\end{lemma}

The next lemma gives an upper bound of a partition's possible main diagonal hook lengths in $\mathcal{DS}(t)$.
 
 \begin{lemma}\label{lem1}
 Let $t$ be a positive integer and $\lambda$ be a partition in $\mathcal{DS}(t)$. If $h \in MD(\lambda)$, then $h \leq 2t-1$.
 \end{lemma}

\pf  Let $\lambda $ be a partition in $\mathcal{DS}(t)$. 
Assume to the contrary that there exists an element $h\in MD(\lambda)$ such that $h > 2t-1$. 
From the definition of $(t,t+1)$-core partitions, we can see that  $h \nequal 2t$ or $2t+2$.  If $h=2t+1 >2t$, then $1 \in MD(\lambda)$ by (\rmnum{1}) of Theorem \ref{pro}. 
However, by (\rmnum{2}) of Theorem \ref{pro},  this is impossible since $1+(2t+1)=2(t+1)$.  Moreover, if $h\geq 2t+3 $, then $h-2t, h-(2t+2) \in MD(\lambda)$ by (\rmnum{1}) of Theorem \ref{pro}. Recall the definition of $d(\lambda)$,   if $h-2t$ is $s$-th main diagonal hook length of $\lambda$, where $1 \leq s \leq |MD(\lambda)|-1$, then we have $d_s=(h-2t)- (h-(2t+2))=2$, contradicting the assumption that  $\lambda \in \mathcal{DS}(t)$. Hence we conclude that $h\leq 2t-1$ for all $h\in MD(\lambda)$ and this completes the proof. \qed

For a positive integer $t$, denote by $\mathrm{Q}_t$  the set $\{1,3,5,7,\ldots,2t-1\} \backslash \{t,t+1\}$.  It is easily checked that $\mathrm{Q}_t$ is a set of distinct positive odd integers. For example, when $t=12$, we have $\mathrm{Q}_{12}=\{1,3,5,7,9,11,15,17,19,21,23\}.$

Combining Theorem \ref{pro} and Lemma \ref{lem1}, we immediately obtain the following characterization of $MD(\lambda)$ for any partition $\lambda$ in $\mathcal{DS}(t)$.

\begin{theorem}\label{ds}
Let $t$ be a positive integer. Then a self-conjugate partition $\lambda$ is a partition in $\mathcal{DS}(t)$ if and only if  $MD(\lambda)$ satisfies the following properties.
\begin{itemize}
\item[\upshape (\rmnum{1}) ] $MD(\lambda) \subseteq \mathrm{Q}_t$; 
\item[\upshape (\rmnum{2}) ]  If $h_1,h_2\in MD(\lambda)$ and $h_1 > h_2$, then $h_1-h_2 \geq4$;
\item[\upshape (\rmnum{3}) ]  If $h_1,h_2\in MD(\lambda)$, then $h_1+h_2 \ne  2t$ and $h_1+h_2 \ne 2t+2$.
\end{itemize}
\end{theorem}

 The following lemma will be essential in proving Theorem \ref{gen1}.
 
\begin{lemma}\label{lrec1}
	For $t\geq3$, we have  the following recurrence relations:
	\begin{align}\label{eq1}
		&f(t)=f(t-1)+f(t-3);\\ \label{eq2} 
		&g(t)=g(t-1)+g(t-3)+f(t-3);\\ \label{eq3}
		&h(t)=2tg(t-1)+(2t-2)g(t-3)-h(t-1)-h(t-3)+f(t-3),
	\end{align}
 where the initial values $g(0)=h(0)=g(1)=h(1)=0$, $f(0)=f(1)=g(2)=h(2)=1$, and $f(2)=2$.

\end{lemma}
 \pf Let $\lambda$ be a partition in $ \mathcal{DS}(t)$ with $MD(\lambda)=\{x_1,x_2,\dots,x_m\}$, where  $x_i>x_{i+1}$ for $1\leq i \leq m-1$. By the definition of $f(t)$, $g(t)$, and $h(t)$, and Lemma \ref{sclar}, we obtain that 
\begin{equation}\label{refdf}
    f(t)=\sum_{\lambda \in \mathcal{DS}(t)}1, \,\,\,\,\,
g(t)=\sum_{\lambda \in \mathcal{DS}(t)}|MD(\lambda)|, \,\,\,\, \, \mbox{and} \,\,\,
h(t)=\sum_{\substack{\lambda \in \mathcal{DS}(t) \\ x\in MD(\lambda)}}x.
\end{equation}
 Now we proceed to distinguish two cases as follows.
 
 \noindent {\bf  Case 1.}  If $1\notin MD(\lambda)$. 
 Let $\Omega(t)$ denote the set of  partitions $\lambda \in \mathcal{DS}(t)$  with $MD(\lambda)\subseteq\mathrm{Q}_{t} \backslash \{1\} $. 
 Now we describe  a map $\psi_{t}: \Omega(t)\rightarrow \mathcal{DS}(t-1) $ as follows. 
   If $\lambda = \emptyset$, then $\psi_{t}(\lambda)=\emptyset$; if $\lambda \ne \emptyset$, then we define $\psi_{t}(\lambda)$ to be the partitions $\mu$ with $MD(\mu)=\{y_1,y_2,\dots,y_m\}$, where
 $y_i=2t-x_{m+1-i}$ for $1\leq i \leq m$.

\ \noindent {\bf  Claim 1.} $ MD(\mu) \subseteq \mathrm{Q}_{t-1}$.\\
 If not, then there exists an element $p \in MD(\mu)$ such that $p \notin \mathrm{Q}_{t-1}$. By the map $\psi_{t}$, we have $2t-p \in MD(\lambda) \subseteq \mathrm{Q}_{t} \backslash \{1\} $. This means that  $2t-p$ is odd,  $3\leq 2t-p \leq 2t-1$, and $2t-p \ne t,t+1$. Then, it is easily seen that $p$ is odd, $1\leq p \leq 2t-3$, and $p \ne t-1,t$. From the definition of  $\mathrm{Q}_{t-1}$, we can see that $p \in \mathrm{Q}_{t-1}$, a contradiction. Hence, we have $ MD(\mu) \subseteq \mathrm{Q}_{t-1}$.

\ \noindent {\bf  Claim 2.} $d(\mu) \geq 4$.\\
If not, then there must exists  elements  $y_i, y_j \in MD(\mu)$ ($1 \leq i < j \leq m$)  such that $y_i>y_j$ and $y_i-y_j< 4$. By the map $\psi_{t}$, we have $2t-x_{m+1-i}>2t-x_{m+1-j}$ and $(2t-x_{m+1-i})-(2t-x_{m+1-j})<4$, that is, $x_{m+1-j}>x_{m+1-i}$ and $x_{m+1-j}-x_{m+1-i}<4$.  This contradicts the definition of the partition $\lambda$ in $\mathcal{DS}(t)$. Hence, we have $y_i-y_j\geq 4 $ for all $1 \leq i < j \leq m$, namely, $d(\mu) \geq 4$.
 
\ \noindent {\bf  Claim 3.} $y_i+y_j \ne 2t-2 $ and  $y_i+y_j\ne 2t $ for any $i, j \in [1,m]$ and $i\ne j$.\\
 If not,  there must exists elements $y_i, y_j\in MD(\lambda)$ such that   $y_i+y_j = 2t-2$ or $y_i+y_j = 2t$. Suppose that $y_i+y_j= 2t-2$. By the map $\psi_{t}$, we have $(2t-x_{m+1-i})+(2t-x_{m+1-j})=2t-2$, that is,  $x_{m+1-i}+x_{m+1-j}=2t+2$, which yields a contradiction with the definition of $\lambda$. The proof for the case $y_i+y_j \ne 2t$ is similar and omitted. Hence the claim is proved.

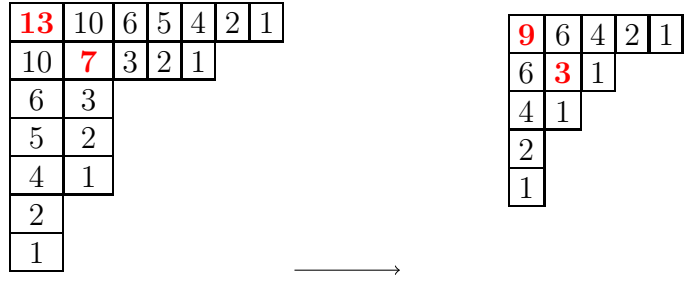
\begin{figure} [h]
 	\centering
 	{$\begin{array}[b]{*{7}c}\cline{1-7}
 			\lr {{\bf\red{13}}}&\lr{10}&\lr{6}&\lr{5}&\lr{4}&\lr{2}&\lr{1}\\\cline{1-7}
 			\lr{10}&\lr{{{\bf\red 7}}}&\lr{3}&\lr{2}&\lr{1}\\\cline{1-5}
 			\lr{6}&\lr{3}\\\cline{1-2}
 			\lr{5}&\lr{2}\\\cline{1-2}
 			\lr{4}&\lr{1}\\\cline{1-2}
 			\lr{2}\\\cline{1-1}
 			\lr{1}\\\cline{1-1}
 		\end{array}$}
 	\begin{tikzpicture}
 		\draw[->] (4.2,0)--(5.6,0);\put(44,1.3){$\psi_{8}(\lambda)$};
 	\end{tikzpicture}
 	\hskip 1.3cm
 	$\begin{array}[b]{*{5}c}\cline{1-5}
 		\lr{{{\bf\red 9}}}&\lr{6}&\lr{4}&\lr{2}&\lr{1}\\\cline{1-5}
 		\lr{6}&\lr{{{\bf\red 3}}}&\lr{1}\\\cline{1-3}
 		\lr{4}&\lr{1}\\\cline{1-2}
 		\lr{2}\\\cline{1-1}
 		\lr{1}\\\cline{1-1}
 		&\\
 		&\\
 	\end{array}$
 	\caption{An example of  Case 1.}\label{fig:md1}
 \end{figure}

For example, let $t=8$ and  $\lambda =(7,5,2,2,2,1) \in \mathcal{DS}(8)$. Then $MD(\lambda)=\{13,7\}$. 
 By applying the map $\psi_{8}$ to $\lambda$, we obtain $MD(\psi_{8}(\lambda))=\{9,3\}$ as shown in Figure \ref{fig:md1}, that is,  $\psi_{8}(\lambda)=(5,3,2,1,1) \in \mathcal{DS}(7)$.
 
 Thus, we conclude that $MD(\mu)$ satisfies the three properties in Theorem \ref{ds}. This means that $\mu \in \mathcal{DS}(t-1)$, that is, the map $\psi_{t}$ is well-defined.

 Conversely, give any $\mu \in \mathcal{DS}(t-1)$, we can recover the partition $\lambda \in \Omega(t)$ as follows.  If $\mu = \emptyset$, then $\lambda=\emptyset \in \Omega(t)$. If $\mu \ne \emptyset$, then we suppose $MD(\mu)=\{y_1,y_2,\dots,y_m\}$. Define the partition $\lambda$ with $MD(\lambda)=\{x_1,x_2,\dots, x_m\}$, then we can recover a partition $\lambda \in \Omega(t)$ by reversing the procedure in Case 1, that is, $x_{i}=2t-y_{m+1-i}$ for all $1\leq i \leq m$. So the construction of the map $\psi_{t}$ is reversible, hence it is a bijection. Moreover, it is apparent that $MD(\lambda)$ satisfies the three properties in Theorem \ref{ds}, which implies that  $\lambda \in \mathcal{DS}(t)$ with $1 \notin MD(\lambda)$.

 \noindent {\bf  Case 2.} If  $1\in MD(\lambda)$. From (\rmnum{2}) and (\rmnum{3}) of Theorem \ref{ds}, then we derive that $3,2t-1\notin MD(\lambda)$ since $3-1=2<4$ and $(2t-1)+1=2t$.   Let $\Omega'(t)$ denote the set of  partitions $\lambda'$, where $\lambda'$ is the partition obtained from  $\lambda$ by deleting the main diagonal hook length $1$, namely, $MD(\lambda')=MD(\lambda) \backslash\{1\}=\{x_1,x_2,\dots,x_{m-1}\} $. 
 Then we have $MD(\lambda')=MD(\lambda) \backslash\{1\} \subseteq \mathrm{Q}_{t} \backslash\{1,3,2t-1\}$.  
 Now we describe a map $\phi_{t}: \Omega'(t) \rightarrow \mathcal{DS}(t-3) $ as follows. 
 If $\lambda' = \emptyset$, that is, $\lambda=(1)$, then $\phi_{t}(\lambda')=\emptyset$; if $\lambda' \ne \emptyset$, then we define $\phi_{t}(\lambda')$ to be the partition $\mu$ with  $MD(\mu)=\{y_1,y_2,\dots,y_{m-1}\}$, where
 $y_i=2t-2-x_{m-i}$ for $1\leq i \leq m-1$.

\ \noindent{\bf Claim 1.} $MD(\mu) \subseteq \mathrm{Q}_{t-3}$.\\
If not, then there exists an element $p \in MD(\mu)$ such that $p \notin \mathrm{Q}_{t-3}$.  By the map $\phi_{t}$, we have $2t-2-p \in MD(\lambda') \subseteq \mathrm{Q}_{t} \backslash \{1,3,2t-1\} $. This means that  $2t-2-p$ is odd,  $5\leq 2t-2-p \leq 2t-3$, and $2t-2-p \ne t,t+1$. Then, it is easily seen that $p$ is odd, $1\leq p \leq 2t-7$, and $p \ne t-3,t-2$. From the definition of $\mathrm{Q}_{t-3}$, we can see that $p \in \mathrm{Q}_{t-3}$, a contradiction. Hence, we have $ MD(\mu) \subseteq \mathrm{Q}_{t-3}$.

\ \noindent{\bf Claim 2.} $d(\mu)\geq4$.\\
If not, then there must exists elements $y_i, y_j \in MD(\mu)$  ($1 \leq i < j \leq m-1$)   such that $y_i>y_j$ and $y_i-y_j< 4$. By the map $\phi_{t}$, we have $2t-2-x_{m-i}>2t-2-x_{m-j}$ and $(2t-2-x_{m-i})-(2t-2-x_{m-j})<4$, that is, $x_{m-j}>x_{m-i}$ and $x_{m-j}-x_{m-i}<4$. Then we have  $x_{m-j}, x_{m-i} \in MD(\lambda)$ since $x_{m-j}, x_{m-i} \in MD(\lambda')$ and $MD(\lambda')=MD(\lambda) \backslash\{1\} $. This contradicts the definition of the partition $\lambda$ in $\mathcal{DS}(t)$. Hence, we have $y_i-y_j\geq 4 $ for all $1 \leq i < j \leq m-1$, that is, $d(\mu) \geq 4$.

\  \noindent{\bf Claim 3.} $y_i+y_j \ne 2t-6$ and  $y_i+y_j \ne 2t-4 $ for any $i, j \in [1,m-1]$ and $i \ne j$.\\
 If not,  there must exists elements  $y_i, y_j \in MD(\mu)$ such that  $y_i+y_j = 2t-6 $ or $y_i+y_j =2t-4$.  Suppose that $y_i+y_j =2t-6$. By the map $\phi_{t}$,  we have $(2t-2-x_{m-i})+(2t-2-x_{m-j})=2t-6$, that is, $x_{m-i}+x_{m-j}=2t+2$.  Then we have $x_{m-j}, x_{m-i} \in MD(\lambda)$ since $x_{m-j}, x_{m-i} \in MD(\lambda')$ and $MD(\lambda')=MD(\lambda) \backslash\{1\} $. Thus, this yields a contradiction with the definition of $\lambda$. The proof for the case $y_i+y_j \ne 2t-4$ is similar and omitted. Hence the claim is proved.

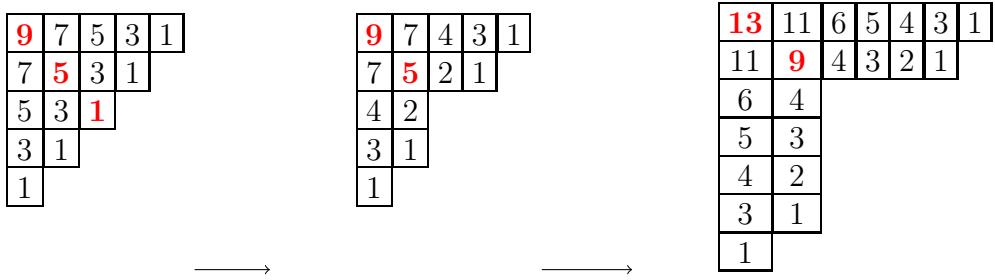
\begin{figure} [h]
 	\centering
 	{$\begin{array}[b]{*{5}c}\cline{1-5}
 			\lr{{\bf\red 9}}&\lr{7}&\lr{5}&\lr{3}&\lr{1}\\\cline{1-5}
 			\lr{7}&\lr{{{\bf\red 5}}}&\lr{3}&\lr{1}\\\cline{1-4}
 			\lr{5}&\lr{3}&\lr{{{\bf \red 1}}}\\\cline{1-3}
 			\lr{3}&\lr{1}\\\cline{1-2}
 			\lr{1}\\\cline{1-1}
 			&\\
 			&\\
 		\end{array}$}
   \begin{tikzpicture}
 		\draw[->] (4.5,0)--(5.5,0);\put(44,1.3){$\lambda \rightarrow \lambda' $};
 	\end{tikzpicture}
 	\hskip 1cm
        $\begin{array}[b]{*{5}c}\cline{1-5}
 			\lr{{\bf\red 9}}&\lr{7}&\lr{4}&\lr{3}&\lr{1}\\\cline{1-5}
 			\lr{7}&\lr{{{\bf\red 5}}}&\lr{2}&\lr{1}\\\cline{1-4}
 			\lr{4}&\lr{2}\\\cline{1-2}
 			\lr{3}&\lr{1}\\\cline{1-2}
 			\lr{1}\\\cline{1-1}
 			&\\
 			&\\
 		\end{array}$
   \begin{tikzpicture}
 		\draw[->] (4.4,0)--(5.6,0);\put(44,1.3){$\phi_{10}(\lambda')$};
 	\end{tikzpicture}
             \hskip 1cm
 	$\begin{array}[b]{*{7}c}\cline{1-7}
 		\lr {{\bf\red{13}}}&\lr{11}&\lr{6}&\lr{5}&\lr{4}&\lr{3}&\lr{1}\\\cline{1-7}
 		\lr{11}&\lr{{{\bf\red 9}}}&\lr{4}&\lr{3}&\lr{2}&\lr{1}\\\cline{1-6}
 		\lr{6}&\lr{4}\\\cline{1-2}
 		\lr{5}&\lr{3}\\\cline{1-2}
 		\lr{4}&\lr{2}\\\cline{1-2}
 		\lr{3}&\lr{1}\\\cline{1-2}
 		\lr{1}\\\cline{1-1}
 	\end{array}$
 	
 	\caption{An example of  Case 2.}\label{fig:md2}
 \end{figure}

For example, let $t=10$ and  $\lambda=(5,4,3,2,1) \in \mathcal{DS}(10)$.  Then $MD(\lambda)=\{9,5,1\}$. Thus, we have $\lambda'=(5,4,2,2,1)$ and $MD(\lambda')=\{9,5\}$. By applying the map $\phi_{10}$ to $\lambda'$, we obtain $MD(\phi_{10}(\lambda'))=\{13,9\}$ as shown in Figure \ref{fig:md2}, that is,  $\phi_{10}(\lambda')=(7,6,2,2,2,2,1) \in \mathcal{DS}(7)$.
 
 Thus, we conclude that $MD(\mu)$ satisfies the three properties in Theorem \ref{ds}. This implies that $\mu \in \mathcal{DS}(t-3)$, that is, the map $\phi_{t}$ is well-defined.
 
Conversely, give any $\mu \in \mathcal{DS}(t-3)$, we can recover the partition $\lambda \in \Omega'(t)$ as follows.  If $\mu = \emptyset$, then $\lambda=\emptyset \in \Omega'(t)$. If $\mu \ne \emptyset$, then we suppose $MD(\mu)=\{y_1,y_2,\dots,y_m\}$. Define the partition $\lambda$ with $MD(\lambda)=\{x_1,x_2,\dots, x_m\}$, then we can recover a partition $\lambda \in \Omega'(t)$ by reversing the procedure in Case 2, that is, $x_{i}=2t-2-y_{m+1-i}$ for all
$1 \leq i \leq m$.  So the construction of the map $\phi_{t}$ is reversible, hence it is a bijection. Furthermore, let $\tau$ be a partition with its corresponding $MD(\tau) = MD(\lambda) \bigcup \{1\}= \{x_1,x_2,\dots, x_m, 1\}$.  One can easily check that $MD(\tau)$ satisfies the three properties in Theorem \ref{ds}, which implies that $\tau \in \mathcal{DS}(t)$  with $1 \in MD(\tau)$.

In the following, we aim to prove (\ref{eq1}), (\ref{eq2}) and (\ref{eq3}). 
% We shall only give a proof of (\ref{eq3}), relations (\ref{eq1}) and (\ref{eq2}) can be verified in the same manner. 
Combining the above two cases and (\ref{refdf}),  we immediately obtain the following results:
 \begin{align*}
 f(t)&=\sum_{\lambda \in \mathcal{DS}(t)}1=\sum_{\substack{\lambda \in \mathcal{DS}(t) \\ 1 \in MD(\lambda)}  } 1+\sum_{\substack{\lambda \in \mathcal{DS}(t) \\ 1\notin MD(\lambda) }} 1=\sum_{\lambda \in \mathcal{DS}(t-3)}1+\sum_{\lambda \in \mathcal{DS}(t-1)}1=f(t-3)+f(t-1)
 \end{align*}
 %as desired, 
 and
\begin{align*}
 g(t)&=\sum_{\lambda \in \mathcal{DS}(t)}|MD(\lambda)|=\sum_{\substack{\lambda \in \mathcal{DS}(t) \\ 1 \in MD(\lambda)}  } |MD(\lambda)|+\sum_{\substack{\lambda \in \mathcal{DS}(t) \\ 1\notin MD(\lambda) }} |MD(\lambda)|
 \\&
 =\sum_{\lambda \in \Omega'(t)}(|MD(\lambda)|+1)+\sum_{\lambda \in \Omega(t)}|MD(\lambda)|
 % &=\sum_{\lambda \in \mathcal{DS}(t-3)}(|MD(\lambda)|+1)+\sum_{\lambda \in \mathcal{DS}(t-1)}|MD(\lambda)|\\
     =g(t-3)+f(t-3)+g(t-1)
  \end{align*} 
as desired. Similarly, it is routine to check that

\begin{align*}
h(t)&=\sum_{\lambda \in \mathcal{DS}(t)} |\lambda|=\sum_{\substack{\lambda \in \mathcal{DS}(t) \\ 1 \in MD(\lambda)}  } |\lambda|
	+\sum_{\substack{\lambda \in \mathcal{DS}(t) \\ 1\notin MD(\lambda) }} |\lambda|= \sum_{\lambda \in \Omega'(t)}(|\lambda|+1)
	+\sum_{\lambda \in \Omega(t)}|\lambda|\\
	&= \sum_{\substack{\lambda \in \mathcal{DS}(t-3)\\ x\in MD(\lambda)}}(2t-2-x)
	+\sum_{\lambda \in \mathcal{DS}(t-3)}1
	+\sum_{\substack{\lambda \in \mathcal{DS}(t-1)\\ x\in MD(\lambda)}}(2t-x) \\
	&=2tg(t-1)+(2t-2)g(t-3)-h(t-1)-h(t-3)+f(t-3)
\end{align*}
as desired, and hence the proof is complete. \qed

Now we are ready to prove Theorem \ref{gen1}.

\noindent{\bf{Proof of Theorem \ref{gen1}.}} 
Let
$$
F(x)=\sum_{t\geq0}f(t-1)x^t, \,\,\,\,\,
G(x)=\sum_{t\geq0}g(t-1)x^t,\,\,\,\,\,
\mbox{and} \,\,\,
H(x)=\sum_{t\geq0}h(t-1)x^t
$$ and set $f(-1)=1$, $g(-1)=h(-1)=0$.
From the recurrence relation (\ref{eq1}) and the initial values $f(-1)=f(0)=f(1)=1$ and $f(2)=2$, we obtain that
$$
F(x)=1+xF(x)+x^3F(x).
$$
which implies that  \eqref{egen1} as desired. 

Similarly, from the recurrence relation \eqref{eq2} and the initial values $g(-1)=g(0)=g(1)=0$ and $g(2)=1$, we deduce that
$$
G(x)=xG(x)+x^3G(x)+x^3F(x)
$$
which implies that 
\begin{equation}\label{geq}
	G(x)=\frac{x^3F(x)}{1-x-x^3}.
\end{equation}
Next, substituting (\ref{egen1}) into (\ref{geq}), we get (\ref{egen2}) as desired.

Furthermore,  from the recurrence relation (\ref{eq3}) and the initial values $h(-1)=h(0)=h(1)=0$ and $h(2)=1$, we derive that
$$
H(x)=2x^3G(x)+(2x^2+2x^4)G'(x)-xH(x)-x^3H(x)+x^3F(x),
$$
which implies that
\begin{equation}\label{heq}
	H(x)=\frac{2x^3G(x)+(2x^2+2x^4)G'(x)+x^3F(x)}{1+x+x^3}.
\end{equation}
Substituting (\ref{egen1})  and (\ref{egen2}) into (\ref{heq}), we obtain (\ref{egen3}) and complete the proof.\qed

\section{Proof of  Theorem  \ref{largest}}
This section is devoted to proving Theorem \ref{largest}. To this end, we shall construct a partition $\lambda$ of $MD(\lambda)$ through the following five lemmas. It turns out that $\lambda$ is the unique self-conjugate $(t,t+1)$-core partitions in $\mathcal{DS}(t)$ which has the largest size.

\begin{lemma}\label{lle1}
For a positive integer $t \geq 3$, let $\lambda$ be a partition in $\mathcal{DS}(t)$ with the largest size. Then $2t-1 \in MD(\lambda)$. Furthermore,  $1,3,2t-3\notin MD(\lambda)$.
\end{lemma}

\pf  Suppose that $\lambda$ is a partition in $\mathcal{DS}(t)$ with the largest size, where $t \geq 3$. 
Let $MD(\lambda)=MD'(\lambda) \bigcup MD''(\lambda)$, where $MD'(\lambda) \subseteq \mathrm{Q}_t \backslash \{1,3,2t-3,2t-1\}$ and $MD''(\lambda) \subseteq \{1,3,2t-3,2t-1\}.$ 
Assume to the contrary that the lemma is not valid, that is,
$MD''(\lambda) \subseteq \{1,3,2t-3\}$.  Since $3+(2t-3)=2t$, by (\rmnum{3}) of Theorem \ref{ds} we obtain that $3$ and $2t-3$ cannot simultaneously belong to $MD''(\lambda)$. Also, $t\geq 3$ implies $2t-3\geq 3$. Therefore, 
$$
\sum_{x\in MD''(\lambda)} x \leq 1+(2t-3).
$$
%If not, we have $3, 2t-3 \in MD''(\lambda)$, which means that $3, 2t-3 \in MD(\lambda)$. However,  by (\rmnum{3}) of Theorem \ref{ds}, this is impossible since $3+(2t-3)=2t$.  By (\rmnum{2}) of Theorem \ref{ds}, we obtain that $1$ and $3$ cannot simultaneously belong to $MD''(\lambda)$ since $3-1=2 < 4$. 
%Since $t \geq 3$, we have $MD''(\lambda)=\{1, 2t-3\}$ guarantee that $\lambda$ is of maximum size. }
Let $\mu$ be a self- conjugate partition  with its corresponding $MD(\mu)=MD'(\lambda) \bigcup \{2t-1\}$.  Since $t\geq 3$, we have $2t-1\neq 1$ or $3$. Thus, it is easy to check that $MD(\mu)$ satisfies the three properties in Theorem \ref{ds}. This means that $\mu \in \mathcal{DS}(t)$.  By Lemma \ref{sclar},  we have 
$$|\lambda| \leq~ (1+(2t-3))+ \sum_{x\in MD'(\lambda)}x~<~(2t-1) + \sum_{x\in MD'(\lambda)}x~=~|\mu|.$$ 
This contradicts the fact that $\lambda$ is of maximum size, thus we have $2t-1 \in MD(\lambda)$.   Furthermore,  by (\rmnum{2}) and (\rmnum{3}) of Theorem \ref{ds},  we have $1,3,2t-3\notin MD(\lambda)$ since $1+(2t-1)=2t$, $3+(2t-1)=2(t+1)$, and $(2t-1)-(2t-3) = 2 < 4$. This completes the proof.\qed

\begin{figure}[H]
	\begin{center}
		\begin{tikzpicture}[font =\small , scale = 1, line width = 0.7pt]
			\filldraw[fill=white](0,0)circle(0.1);
			\filldraw[fill=white](0.8,0)circle(0.1);
			\filldraw[fill=white](1.6,0)circle(0.1);
			\filldraw[fill=white](2.4,0)circle(0.1);
			\filldraw[fill=white](3.2,0)circle(0.1);
						
			\filldraw[fill=white](-0.4,-1.5)circle(0.1);
			\filldraw[fill=white](0.4,-1.5)circle(0.1);
			\filldraw[fill=white](1.2,-1.5)circle(0.1);
			\filldraw[fill=white](2,-1.5)circle(0.1);
			\filldraw[fill=white](2.8,-1.5)circle(0.1);
			\filldraw[fill=white](3.6,-1.5)circle(0.1);
						
			  \draw[very thick](0,0)--(-0.4,-1.5);
			  \draw[very thick](0,0)--(0.4,-1.5);
			\draw[very thick](0.8,0)--(0.4,-1.5);
			  \draw[very thick](0.8,0)--(1.2,-1.5);
			\draw[very thick](1.6,0)--(1.2,-1.5);
			  \draw[very thick](1.6,0)--(2,-1.5);
			\draw[very thick](2.4,0)--(2,-1.5);
			\draw[very thick](2.4,0)--(2.8,-1.5);
			\draw[very thick](3.2,0)--(2.8,-1.5);
			\draw[very thick](3.2,0)--(3.6,-1.5);

			 \coordinate [label=above:$15$] (x) at (0,0.1);
			 \coordinate [label=above:$17$] (x) at (0.8,0.1);
			 \coordinate [label=above:$19$] (x) at (1.6,0.1);
			 \coordinate [label=above:$21$] (x) at (2.4,0.1);
			 \coordinate [label=above:$23$] (x) at (3.2,0.1);
			 \coordinate [label=above:$11$] (x) at (-0.4,-2.2);
			 \coordinate [label=above:$9$] (x) at (0.4,-2.2);
		   \coordinate [label=above:$7$] (x) at (1.2,-2.2);
	       \coordinate [label=above:$5$] (x) at (2,-2.2);
		   \coordinate [label=above:$3$] (x) at (2.8,-2.2);
              \coordinate [label=above:$1$] (x) at (3.6,-2.2);
		
		\end{tikzpicture}
	\end{center}
	\caption{The  representation of  $\mathrm{Q}_{12}^{l}$ and $\mathrm{Q}_{12}^{r}.$}\label{fig:qlr}
\end{figure}
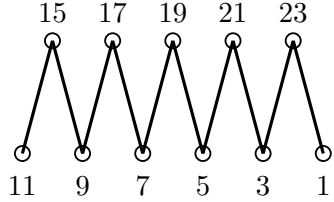

Based on the previous definition, we now divide the set $\mathrm{Q}_{t}$ into two subsets. Define the set $\mathrm{Q}_t^{l}=\mathrm{Q}_{t} \bigcap [t+2,2t-1]$ and $\mathrm{Q}_t^{r}=\mathrm{Q}_{t} \bigcap [1, t-1]$.
It is easily checked that $\mathrm{Q}_t=\mathrm{Q}_t^{l} \bigcup \mathrm{Q}_t^r$. For the example mentioned above, when $t=12$, we have $\mathrm{Q}_{12}^{l}=\{15,17,19,21,23\}$ and $\mathrm{Q}_{12}^{r}=\{1,3,5,7,9,11\}.$ Figure \ref{fig:qlr} lays out the
 representations of $\mathrm{Q}_{12}^{l}$ and $\mathrm{Q}_{12}^{r}.$ There is an edge between two numbers in Figure \ref{fig:qlr} iff their summation is equal to $2t$ or $2t+2$.

In the following, we give some descriptions for $d(\lambda)$ of a partition $\lambda \in \mathcal{DS}(t)$ with the largest size. For partition $\lambda \in \mathcal{DS}(t)$ with the largest size, notice that $t \geq 6$ if the sequence $d(\lambda)$ is not empty.

\begin{lemma}\label{lle2}
For a positive integer $t \geq 6$, let $\lambda$ be a partition in $\mathcal{DS}(t)$ with the largest size. If $h_i, h_{i+1} \in MD(\lambda)$ for $1\leq i \leq |MD(\lambda)|-1$, then $4 \leq d_i \leq 6$.
\end{lemma}

\pf Suppose that $\lambda$ is a partition in $\mathcal{DS}(t)$ with the largest size,  where $t \geq 6$. Assume to the contrary that  $MD(\lambda)$ contain elements $h_j, h_{j+1}$ such that $d_{j}\geq8$, that is, $h_j-h_{j+1}\geq 8$. Choose $j$ to be the
smallest such integer. Then we consider the following cases.

\noindent {\bf  Case 1.} If $h_{j}, h_{j+1}\in\mathrm{Q}_{t}^{l}$.\\
\noindent {\bf  Subcase 1.1.} If $2t+4-h_j \, (\mbox{resp.} \, 2t+6-h_j) \in MD(\lambda)$, then we  obtain a new partition $\lambda'$, whose corresponding $MD(\lambda')$ is  obtained from $MD(\lambda)$ by  replacing  $2t+4-h_j \, (\mbox{resp.}\, 2t+6-h_j) $ by $h_j-4$. \\
\noindent {\bf  Subcase 1.2.} If $2t+4-h_j, 2t+6-h_j \notin MD(\lambda)$, then we  obtain a new partition $\lambda'$, whose corresponding $MD(\lambda')$ is  obtained from $MD(\lambda)$ by adding  $h_j-4$.

For example, let $\lambda \in \mathcal{DS}(16)$ with $MD(\lambda) = \{31,23,19,5\}$. Then $d(\lambda)=(8,4,6)$ and $h_1=31, h_2=23 \in \mathrm{Q}_{16}^{l}$. By the above procedure in Subcase 1.1, we obtain $\lambda' \in \mathcal{DS}(16)$ with $MD(\lambda')=\{31,27,23,19\}$.

\noindent {\bf  Case 2.} If $h_j\in\mathrm{Q}_{t}^{l}$ and $h_{j+1}\in\mathrm{Q}_{t}^{r}$, then we have $2t-h_{j+1} \notin MD(\lambda)$.\\
\noindent {\bf  Subcase 2.1.}  When $t$ is even, then $t+3$ is the smallest element of  $\mathrm{Q}_{t}^{l}$. Thus, we further distinguish as follows.
\begin{itemize}
\item[(1)] If $h_j=t+3$, then we obtain a new partition $\lambda'$, whose corresponding $MD(\lambda')$ is  obtained from $MD(\lambda)$ by replacing each $x$ by $x+2$ for $x\in MD(\lambda)$ and $h_j\leq x \leq 2t-2-h_{j+1}$,  and replacing  $h_{j+1}$ by $t-1$.
\item[(2)] If $h_j\ne t+3$ and $h_{j+1} \ne t-1, t-3$, then we obtain a new partition $\lambda'$, whose corresponding $MD(\lambda')$ is  obtained from $MD(\lambda)$ by adding $t-1$.
\item[(3)] If $h_j\ne t+3$ and $h_{j+1} = t-1\,\, \mbox{or} \,\, t-3$, then we have $t+5 \notin MD(\lambda)$. Thus, we can obtain a new partition $\lambda'$, whose corresponding $MD(\lambda')$ is  obtained from $MD(\lambda)$ by replacing $h_{j+1}$ by $t+3$.
\end{itemize}

For example, let $\lambda \in \mathcal{DS}(18)$ with $MD(\lambda) = \{35,31,25,21,9\}$. Then $d(\lambda)=(4,6,4,12)$ and $h_4=21 \in \mathrm{Q}_{18}^{l}, \, h_5=9 \in \mathrm{Q}_{18}^{r}$. By the above procedure in Subcase 2.1(1), we obtain $\lambda' \in \mathcal{DS}(18)$ with $MD(\lambda')=\{35,31,27,23,17\}$.

\noindent {\bf  Subcase 2.2.} When $t$ is odd, then $t+2$ is the smallest element of  $\mathrm{Q}_{t}^{l}$. Thus, we further distinguish as follows.
\begin{itemize}
\item[(1)] If $h_j=t+2$, then we obtain a new partition $\lambda'$, whose corresponding $MD(\lambda')$ is obtained from $MD(\lambda)$ by replacing each   $x$ by $x+2$ for $x\in MD(\lambda)$ and $h_j+4\leq x \leq 2t-2-h_{j+1}$,  and replacing  $h_{j+1}$ by $t-4$.
\item[(2)] If $h_j\ne t+2$ and $h_{j} < 2t-h_{j+1}$, then we obtain a new partition $\lambda'$, whose corresponding $MD(\lambda')$ is obtained from $MD(\lambda)$ by replacing $x$ by $x+2$ for $h_j\leq x\leq 2t-2-h_{j+1}$ and $x\in MD(\lambda)$, and replacing $h_{j+1}$ by $t+2$.
\item[(3)] If $h_j\ne t+2$ and $h_{j} > 2t-h_{j+1}$, then we obtain a new partition $\lambda'$, whose corresponding $MD(\lambda')$ is obtained from $MD(\lambda)$ by replacing $h_{j+1}$ by $t+2$.
\end{itemize}

For example, let $\lambda \in \mathcal{DS}(15)$ with $MD(\lambda) = \{29,25,19,9\}$. Then we have $d(\lambda)=(4,6,10)$, $h_3=19 \in \mathrm{Q}_{15}^{l}, \, h_4=9 \in \mathrm{Q}_{15}^{r}$, and $h_3 =19 < 2\times 15-h_{4}=21$. By the above procedure in Subcase 2.2(2), we obtain $\lambda'$ with $MD(\lambda')=\{29,25,21,17\}$.

\noindent {\bf  Case 3.} If $h_j, h_{j+1}\in\mathrm{Q}_{t}^{r}$, then we obtain that $2t-h_{j+1}, 2t+2-h_{j+1} \notin MD(\lambda)$. Also, by Lemma \ref{lle1} and $h_j-h_{j+1}\geq 8$, we know that $5\leq h_{j+1}\leq t-9$.
By choice of $j$, we know that $4 \leq d_i \leq 6$ for all $i < j$. This implies that  $2t+4-h_{j+1}, 2t-2-h_{j+1} \in MD(\lambda)$. Furthermore, we have $h_{j+1} +4 \notin MD(\lambda)$ since $2t+4-h_{j+1} \in MD(\lambda)$.
Thus, we can obtain a new partition $\lambda'$, whose corresponding $MD(\lambda')$ is obtained from $MD(\lambda)$ by replacing $2t-2-h_{j+1}$ by $2t-h_{j+1}$, and replacing $h_{j+1}$ by $h_{j+1}+4$.

For example, let $\lambda \in \mathcal{DS}(19)$ with $MD(\lambda) = \{37,31,27,21,15,5\}$. Then we have $d(\lambda)=(6,4,6,6,10)$ and $h_5=15, h_6=5 \in \mathrm{Q}_{19}^{r}$. By the above procedure in Case 3, we obtain $\lambda'  \in \mathcal{DS}(19)$ with $MD(\lambda')=\{37,33,27,21,15,9\}$.

 By the above analysis, it is easy to check that $MD(\lambda')$ satisfies the three properties in Theorem \ref{ds} and has a larger sum of the elements as $MD(\lambda)$. From Lemma \ref{sclar},  we obtain that the partition size corresponding to $MD(\lambda')$ is larger than that of $\lambda$, which contradicts the assumption that $\lambda$ has the largest size. This completes the proof.\qed

\begin{lemma}\label{llea}
For a positive integer  $t\geq 6$, let $\lambda$ be a partition in $\mathcal{DS}(t)$ with the largest size. If $h_i, h_{i+1} \in MD(\lambda)\bigcap \mathrm{Q}_{t}^{l}$ and $d_i=6$ for $1\leq i \leq |MD(\lambda)|-1$, then $2t+4-h_{i} \in MD(\lambda)$. 
\end{lemma}
\pf Assume to the contrary the lemma is not valid, that is, there exist elements $h_i,h_{i+1} \in MD(\lambda)\bigcap \mathrm{Q}_{t}^{l}$ such that $d_i=6$ (which means that $h_i=h_{i+1}+6$) and $2t+4-h_{i}=2t-2-h_{i+1} \notin MD(\lambda)$. Therefore, $2t+6-h_{i}=2t-h_{i+1}\notin MD(\lambda)$ and $2t+2-h_{i}\notin MD(\lambda)$. Then we can obtain a new partition $\lambda'$, whose corresponding $MD(\lambda')$ is obtained from $MD(\lambda)$ by adding $2t+4-h_{i}$. 
%An illustration of this process in a generic small example is shown in Figure \ref{fig:lle2}. 
It is easily seen that $MD(\lambda')$ satisfies the three properties in Theorem \ref{ds} and has a larger sum of the elements as $MD(\lambda)$. By Lemma \ref{sclar}, the partition size corresponding to $MD(\lambda')$ is larger than that of $\lambda$. This contradicts that $\lambda$ has the largest size, completing the proof.\qed

\begin{lemma}\label{lle3}
For a positive integer  $t\geq 6$, let $\lambda$ be a partition in $\mathcal{DS}(t)$ with the largest size. If $h_i, h_{i+1} \in MD(\lambda)$ and $d_i=6$ for $1\leq i \leq |MD(\lambda)|-1$, then $d_j=6$ for all $j>i$.
\end{lemma}
\pf Assume that $\lambda$ is a partition in $\mathcal{DS}(t)$ with the largest size,  where $t \geq 6$. Let $h_i, h_{i+1} \in MD(\lambda)$ with  $d_i=6$. By Lemma \ref{lle2}, we proceed to show that there does not exist   $q>i$ such that $h_q, h_{q+1} \in MD(\lambda)$ with $d_q=4$. If not, let $q$ be the smallest such integer such that $h_q,h_{q+1} \in MD(\lambda)$ with $d_q=4$.  By the choice of $q$, we derive that $d_{q-1}=6$ and $h_{q-1}, h_{q} \in \mathrm{Q}_{t}^{l}$.
From Lemma \ref{llea}, we have $2t+4-h_{q-1}=2t-2-h_q \in MD(\lambda)$. Define $\lambda'$ to be the partition corresponding $MD(\lambda')$, where $MD(\lambda')$ is obtained from $MD(\lambda)$ by replacing $h_q$ by $h_q+2$, and replacing $2t-2-h_q$ by $2t+2-h_q$. One can easily derive that  $MD(\lambda')$ satisfies the three properties in Theorem \ref{ds} and  
has a larger sum of the elements as $MD(\lambda)$. Then we have $|\lambda'|>|\lambda|$, a contradiction with the fact that $\lambda$ is of maximum size.
Hence, there does not exist $q>i$ such that $h_q, h_{q+1} \in MD(\lambda)$ with $d_q=4$. This completes the proof.\qed

\begin{lemma}\label{lle4}
For a positive integer $t\geq 6$, let $\lambda$ be a partition in $\mathcal{DS}(t)$ with the largest size and let $d(\lambda)=(d_1,d_2,\dots,d_{k-1})$, where $k \geq 2$. If there exists a positive integer $1\leq \ell \leq k$ such that $d_{\ell-1}=4$ and $d_{\ell}=6$ (Here $\ell=1$ means that $d_1=d_2=\dots=d_{k-1}=6$; and $\ell=k$ means that $d_1=d_2=\dots=d_{k-1}=4$), then 
\begin{equation}
	\lfloor \frac{2t+1}{12} \rfloor  \leq \ell \leq \lfloor \frac{2t+1}{12} \rfloor +2.
\end{equation}
\end{lemma}
\pf Suppose that $\lambda$ is a partition in $\mathcal{DS}(t)$ with the largest size and  its corresponding $MD(\lambda) = \{h_1, h_2, \dots, h_k\}$, where $t \geq 6$ and $k \geq 2$. Let $d(\lambda)=(d_1,d_2,\dots,d_{k-1})$ where $d_{\ell-1}=4$ and $d_{\ell}=6$ for some  $1 \leq \ell \leq k$. 
By Lemma \ref{lle2} and Lemma \ref{lle3}, we have $d_1=\dots=d_{\ell-1}=4$ and $d_{\ell}=\dots=d_{k-1}=6$. It can be deduced that $h_{i}=2t+3-4i$ for $1\leq i\leq \ell$  and $h_i=2t+3-4\ell-6(i-\ell)$ for $\ell+1 \leq i\leq k$.
Now we first show that $\ell \geq 	\lfloor \frac{2t+1}{12} \rfloor  $. For a positive integer $1 \leq q \leq \ell$,  we assert that $d_{q-1}=4$ if $h_{q} \geq \frac{4t+8}{3}$. Otherwise, by Lemma \ref{lle2}, we have $d_{q-1}=6$. Let $q$ be the smallest such integer. 
In view of  Lemma \ref{lle3},  we obtain that $d_{q-1}=\dots =d_{k-1}=6$. By Lemma \ref{llea}, we have  $2t+4-h_q, 2t-2-h_q \in MD(\lambda)$. Then we can obtain a new partition $\mu \in \mathcal{DS}(t)$, where $MD(\mu)$ is obtained from $MD(\lambda)$ by replacing $h_q$ by $h_q+2$, replacing $2t+4-h_q$ by $h_q-2$,
and deleting $2t-2-h_{q}$. Since $\lambda$ has the largest size,  we have $|\lambda|>|\mu|$. By Lemma \ref{sclar}, we deduce that $h_q+(2t+4-h_q)+(2t-2-h_q)>(h_q+2)+(h_q-2),$ namely, $h_{q}<\frac{4t+2}{3}$. This contradicts the fact that $h_{q} \geq \frac{4t+8}{3}$.
Hence, we have $d_{q-1}=4$ if $h_{q} \geq \frac{4t+8}{3}$.  Recall that $d_{\ell-1}=4$ and $d_{\ell}=6$, then we have $\ell \geq q$. It yields that  $h_{\ell}=2t+3-4\ell \leq \frac{4t+8}{3}$.  Then we can deduce that $\ell \geq \lfloor \frac{2t+1}{12} \rfloor$ since $\ell$ is a positive integer and the subscripts of $h$ are sorted in ascending order.

Now we proceed to show that $\ell \leq \lfloor \frac{2t+1}{12} \rfloor +2.$ Assume to the contrary that $\ell \geq \lfloor \frac{2t+1}{12} \rfloor +3.$ This means that there exists some $\lfloor \frac{2t+1}{12}\rfloor < p\leq  \ell$ such that $MD(\lambda)$ contains elements $h_p$, $h_p-4$, and $h_{p}-8$.   Choose $p$ to be the largest such integer. Hence, based on the previous proof, one can easily derive that $h_p<\frac{4t+8}{3}$. Define $\mu$ to be a partition in $\mathcal{DS}(t)$, whose  corresponding $MD(\mu)$ is obtained from $MD(\lambda)$ by replacing $h_p$ by $h_p-2$, replacing $h_p-4$ by $2t+6-h_p$, and adding $2t-h_p$. Notice that $|\lambda|>|\mu|$ since $\lambda$ has the largest size. Moreover, by Lemma \ref{sclar}, we deduce that
$h_p+(h_p-4)>(h_p-2)+(2t+6-h_p)+(2t-h_p).$ Hence we have  $h_p>\frac{4t+8}{3}$, a contradiction. Then we conclude that $\ell \leq \lfloor \frac{2t+1}{12} \rfloor +2$ as desired,  completing the proof.\qed

Now we are ready for the proof of Theorem \ref{largest}.

\noindent{\bf{Proof of Theorem \ref{largest}.}} By the above analysis, it is straightforward to provide the characterization for the set of the main diagonal hook lengths of $\lambda \in \mathcal{DS}(t)$ with the largest size, where $t \geq 3$. Recall that $MD(\lambda)=\{h_1,h_2,\dots,h_k\}$ with $h_i>h_{i+1}$ for $ 1 \leq i \leq k-1$.  When $k=1$, that is, $3 \leq t \leq 5$, by Lemma \ref{lle1}, it is easy to check that $MD(\lambda)= \{2t-1\}$. Thus, by Lemma \ref{sclar},  we have $|\lambda|=2t-1$. When $k \geq 2$, that is, $t \geq 6$, let $d(\lambda)=(d_1,d_2,\dots,d_{k-1})$, where $d_{\ell-1}=4$ and $d_{\ell}=6$ for some $1 \leq  \ell \leq  k$. 
%(If $d_1=d_2=\dots=d_{k-1}=4$, then similar to the proof of Lemma \ref{lle4} we can show that $\lambda$ can not have the largest size, a contradiction).
By Lemma \ref{lle4}, we have $d_1=d_2=\dots=d_{\lfloor \frac{2t+1}{12}\rfloor-1}=4$. 
 Then by Lemma \ref{lle1}, we have $h_i=2t+3-4i$ for $1 \leq i\leq \lfloor \frac{2t+1}{12}\rfloor$, namely, $\{2t-1, 2t-5, \dots, 2t+3-4\lfloor \frac{2t+1}{12}\rfloor\} \subseteq MD(\lambda)$. Hence, it is sufficient to determine the value of $\ell$ in Lemma \ref{lle4}.

Without loss of generality,  we assume that $t$ is even and $t \geq 6$. 
% Combining  Lemmas \ref{lle2}, \ref{llea}, \ref{lle3}, and \ref{lle4}, we derive that  $t+5, t-1 \in MD(\lambda)$ when $\ell \ne k$. 
% On the other hand, for a given even positive integer $t$, i
It is not difficult to check that $|\mathrm{Q}_{t}^{l}|=\frac{t}{2}-1$. Recall that $\{2t-1, 2t-5, \dots, 2t+3-4\lfloor \frac{2t+1}{12}\rfloor\} \subseteq MD(\lambda)$. Then we have $2t-3, 2t-7, \dots, 2t+1-4\lfloor \frac{2t+1}{12}\rfloor \notin MD(\lambda)$ by the definition of partitions in $\mathcal{DS}(t)$.
Thus we obtain that
$$|\mathrm{Q}_{t}^{l}\backslash \{2t-1,2t-3,\dots, 2t+1-4\lfloor \frac{2t+1}{12}\rfloor\}|=\frac{t}{2}-1-2\lfloor \frac{2t+1}{12}\rfloor,$$ which is denoted by $\Delta$.
Next, we discuss according to the remainder of $\Delta$ modulo~$3$.

\noindent{\bf Case 1.} When $\Delta=\frac{t}{2}-1-2\lfloor \frac{2t+1}{12}\rfloor \equiv  0 \, (\text{mod}\ 3)$. Combining  Lemmas \ref{lle2}, \ref{llea}, \ref{lle3}, and \ref{lle4}, we derive that  $t+5, t-1 \in MD(\lambda)$ when $\ell \ne k$.
% Recall that $t+5, t-1 \in MD(\lambda)$ when $\ell \ne k$. 
Furthermore,  we deduce that  $2t-3-4\lfloor \frac{2t+1}{12}\rfloor, 2t-9-4\lfloor \frac{2t+1}{12}\rfloor, \dots, t+5, t-1, \dots, 7+4\lfloor \frac{2t+1}{12}\rfloor \in MD(\lambda) $. 
Recall that $2t+3-4\lfloor \frac{2t+1}{12}\rfloor \in MD(\lambda)$. Hence, by Lemma \ref{llea}, it is easily seen that  $1+4\lfloor \frac{2t+1}{12}\rfloor \in MD(\lambda) $ guarantees that $\lambda$ is of maximum size. Then we deduce that
$\ell=\lfloor \frac{2t+1}{12}\rfloor$, namely, $d(\lambda)=(d_1,\dots,d_{\lfloor \frac{2t+1}{12}\rfloor-1},d_{\lfloor \frac{2t+1}{12}\rfloor},\dots, d_{k-1})=(4,\dots,4,6,\dots,6)$.
Thus, we have
$$
MD(\lambda)=\{2t-1,2t-5,\dots,2t+3-4\lfloor \frac{2t+1}{12}\rfloor, 2t-3-4\lfloor \frac{2t+1}{12}\rfloor, \dots, 1+4\lfloor \frac{2t+1}{12}\rfloor \}.
$$
According to Lemma \ref{sclar}, the size of $\lambda$ is 
$$
|\lambda|
=\frac{1}{3}(t-1)\left(t+1-4\lfloor \frac{2t+1}{12} \rfloor \right)+\lfloor \frac{2t+1}{12} \rfloor \left(2t+1-2\lfloor \frac{2t+1}{12} \rfloor \right)
$$
as desire.

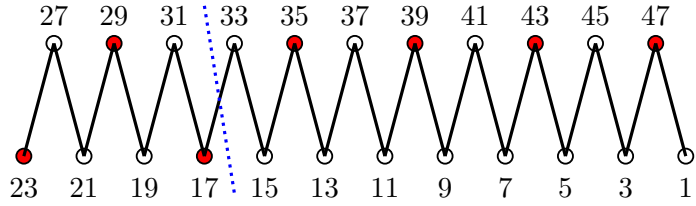
\begin{figure}[H]
	\begin{center}
		\begin{tikzpicture}[font =\small , scale = 1, line width = 0.7pt]
			\filldraw[fill=white](0,0)circle(0.1);
			\filldraw[fill=red](0.8,0)circle(0.1);
			\filldraw[fill=white](1.6,0)circle(0.1);
			\filldraw[fill=white](2.4,0)circle(0.1);
			\filldraw[fill=red](3.2,0)circle(0.1);
			\filldraw[fill=white](4,0)circle(0.1);
			\filldraw[fill=red](4.8,0)circle(0.1);
			\filldraw[fill=white](5.6,0)circle(0.1);
			\filldraw[fill=red](6.4,0)circle(0.1);
			\filldraw[fill=white](7.2,0)circle(0.1);
			\filldraw[fill=red](8,0)circle(0.1);
			
			\filldraw[fill=red](-0.4,-1.5)circle(0.1);
			\filldraw[fill=white](0.4,-1.5)circle(0.1);
			\filldraw[fill=white](1.2,-1.5)circle(0.1);
			\filldraw[fill=red](2,-1.5)circle(0.1);
			\filldraw[fill=white](2.8,-1.5)circle(0.1);
			\filldraw[fill=white](3.6,-1.5)circle(0.1);
			\filldraw[fill=white](4.4,-1.5)circle(0.1);
			\filldraw[fill=white](5.2,-1.5)circle(0.1);
			\filldraw[fill=white](6,-1.5)circle(0.1);
			\filldraw[fill=white](6.8,-1.5)circle(0.1);
			\filldraw[fill=white](7.6,-1.5)circle(0.1);
			\filldraw[fill=white](8.4,-1.5)circle(0.1);
			
			  \draw[very thick](0,0)--(-0.4,-1.5);
			  \draw[very thick](0,0)--(0.4,-1.5);
			\draw[very thick](0.8,0)--(0.4,-1.5);
			  \draw[very thick](0.8,0)--(1.2,-1.5);
			\draw[very thick](1.6,0)--(1.2,-1.5);
			  \draw[very thick](1.6,0)--(2,-1.5);
			  \draw[very thick](2.4,0)--(2,-1.5);
			\draw[very thick](2.4,0)--(2.8,-1.5);
			\draw[very thick](3.2,0)--(2.8,-1.5);
			\draw[very thick](3.2,0)--(3.6,-1.5);
			\draw[very thick](4,0)--(3.6,-1.5);
			\draw[very thick](4,0)--(4.4,-1.5);
			\draw[very thick](4.8,0)--(4.4,-1.5);
		    \draw[very thick](4.8,0)--(5.2,-1.5);
		    \draw[very thick](5.6,0)--(5.2,-1.5);
		    \draw[very thick](5.6,0)--(6,-1.5);
		    \draw[very thick](6.4,0)--(6,-1.5);
		      \draw[very thick](6.4,0)--(6.8,-1.5);
			\draw[very thick](7.2,0)--(6.8,-1.5);
		    \draw[very thick](7.2,0)--(7.6,-1.5);
		    \draw[very thick](8,0)--(7.6,-1.5);
		    \draw[very thick](8,0)--(8.4,-1.5);
		   
               \draw[blue][dotted][very thick](2,0.5)--(2.4,-2);

			 \coordinate [label=above:$27$] (x) at (0,0.1);
			 \coordinate [label=above:$29$] (x) at (0.8,0.1);
			 \coordinate [label=above:$31$] (x) at (1.6,0.1);
			 \coordinate [label=above:$33$] (x) at (2.4,0.1);
			 \coordinate [label=above:$35$] (x) at (3.2,0.1);
			 \coordinate [label=above:$37$] (x) at (4,0.1);
			 \coordinate [label=above:$39$] (x) at (4.8,0.1);
			 \coordinate [label=above:$41$] (x) at (5.6,0.1);
			 \coordinate [label=above:$43$] (x) at (6.4,0.1);
			 \coordinate [label=above:$45$] (x) at (7.2,0.1);
			 \coordinate [label=above:$47$] (x) at (8,0.1);
			 \coordinate [label=above:$23$] (x) at (-0.4,-2.2);
			 \coordinate [label=above:$21$] (x) at (0.4,-2.2);
		   \coordinate [label=above:$19$] (x) at (1.2,-2.2);
	       \coordinate [label=above:$17$] (x) at (2,-2.2);
		   \coordinate [label=above:$15$] (x) at (2.8,-2.2);
              \coordinate [label=above:$13$] (x) at (3.6,-2.2);
		   \coordinate [label=above:$11$] (x) at (4.4,-2.2);
	       \coordinate [label=above:$9$] (x) at (5.2,-2.2);
		   \coordinate [label=above:$7$] (x) at (6,-2.2);
	       \coordinate [label=above:$5$] (x) at (6.8,-2.2);
		   \coordinate [label=above:$3$] (x) at (7.6,-2.2);
	       \coordinate [label=above:$1$] (x) at (8.4,-2.2);
		
		\end{tikzpicture}
	\end{center}
	\caption{The  representation of  $MD(\lambda)=\{47,43,39,35,29,23,17\}$.}\label{fig:lar1}
\end{figure}

For example, let $t=24$ and  $\lambda$ be a partition in $\mathcal{DS}(24)$ with the largest size. Then 
$47,43,39,35 \in MD(\lambda)$. Furthermore, we have $\ell=\lfloor \frac{2\times24+1}{12}\rfloor=4$ and $\Delta=\frac{24}{2}-1-2\times 4=3$. It follows that $29,23, 17\in MD(\lambda)$. Hence we obtain $MD(\lambda)=\{47,43,39,35,29,23,17\}$ as shown in Figure \ref{fig:lar1}. Then $|\lambda|=233$.

\noindent{\bf Case 2.} When $\Delta= \frac{t}{2}-1-2\lfloor \frac{2t+1}{12}\rfloor \equiv  1 \, (\text{mod}\ 3)$.  
% Recall that $t+5, t-1 \in MD(\lambda)$ when $\ell \ne k$.
Combining  Lemmas \ref{lle2}, \ref{llea}, \ref{lle3}, and \ref{lle4}, we derive that  $t+5, t-1 \in MD(\lambda)$ when $\ell \ne k$.
This implies that $\Delta \geq 4$. 
(If $\Delta=1$, we have $t+3 \in MD(\lambda)$ and $\ell =k$, which means that 
$d(\lambda)=(d_1,d_2,\dots,d_{k-1})=(4,4,\dots,4)$.) 
Furthermore,  we deduce that $2t-5-4\lfloor \frac{2t+1}{12}\rfloor, 2t-11-4\lfloor \frac{2t+1}{12}\rfloor, \dots, t+5, t-1, \dots, 9+4\lfloor \frac{2t+1}{12}\rfloor \in MD(\lambda)$. 
Recall that $2t+3-4\lfloor \frac{2t+1}{12}\rfloor \in MD(\lambda)$.  
To guarantee that $\lambda$ is of maximum size, we need to add $2t-1-4\lfloor \frac{2t+1}{12}\rfloor$ to $MD(\lambda)$. Then we have $\ell=\lfloor \frac{2t+1}{12}\rfloor+2$, namely, $d(\lambda)=(d_1,\dots,d_{\lfloor \frac{2t+1}{12}\rfloor+1},d_{\lfloor \frac{2t+1}{12}\rfloor+2},\dots, d_{k-1})=(4,\dots,4,6,\dots,6)$. Thus, we have
$$
MD(\lambda)=\{2t-1, 2t-5, \dots, 2t-5-4\lfloor \frac{2t+1}{12}\rfloor, 2t-11-4\lfloor \frac{2t+1}{12}\rfloor, \dots, 9+4\lfloor \frac{2t+1}{12}\rfloor \}.
$$
% \textcolor{blue}{In addition, if $\Delta=1$, we add $t+3$ to $MD(\lambda)$ guarantees that $\lambda$ is of maximum size.  Then we have $\ell =k$, which means that $d(\lambda)=(d_1, d_2, \dots, d_{k-1})=(4,4,\dots,4)$. Thus, we have 
% \begin{equation}\label{md2}
%    MD(\lambda)=\{2t-1, 2t-5, \dots, 2t+3-4\lfloor \frac{2t+1}{12}\rfloor, t+3\}. 
% \end{equation}
% One can easily check that (\ref{md2}) is a special case of (\ref{md1}). }
According to Lemma \ref{sclar}, the size of $\lambda$ is 
$$
|\lambda|
=\frac{1}{3}(t-1)\left(t-7-4\lfloor \frac{2t+1}{12} \rfloor \right)+\left(\lfloor \frac{2t+1}{12} \rfloor +2\right)\left(2t-3-2\lfloor \frac{2t+1}{12} \rfloor \right)
$$
as desire. 

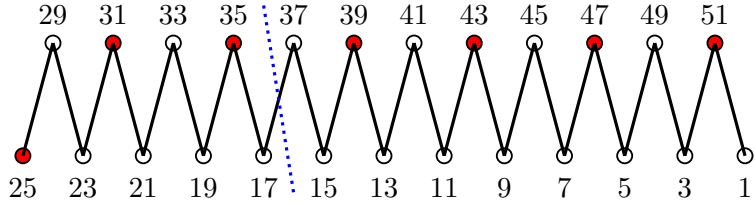
\begin{figure}[H]
	\begin{center}
		\begin{tikzpicture}[font =\small , scale = 1, line width = 0.7pt]
			\filldraw[fill=white](0,0)circle(0.1);
			\filldraw[fill=red](0.8,0)circle(0.1);
			\filldraw[fill=white](1.6,0)circle(0.1);
			\filldraw[fill=red](2.4,0)circle(0.1);
			\filldraw[fill=white](3.2,0)circle(0.1);
			\filldraw[fill=red](4,0)circle(0.1);
			\filldraw[fill=white](4.8,0)circle(0.1);
			\filldraw[fill=red](5.6,0)circle(0.1);
			\filldraw[fill=white](6.4,0)circle(0.1);
			\filldraw[fill=red](7.2,0)circle(0.1);
			\filldraw[fill=white](8,0)circle(0.1);
			\filldraw[fill=red](8.8,0)circle(0.1);
			
			\filldraw[fill=red](-0.4,-1.5)circle(0.1);
			\filldraw[fill=white](0.4,-1.5)circle(0.1);
			\filldraw[fill=white](1.2,-1.5)circle(0.1);
			\filldraw[fill=white](2,-1.5)circle(0.1);
			\filldraw[fill=white](2.8,-1.5)circle(0.1);
			\filldraw[fill=white](3.6,-1.5)circle(0.1);
			\filldraw[fill=white](4.4,-1.5)circle(0.1);
			\filldraw[fill=white](5.2,-1.5)circle(0.1);
			\filldraw[fill=white](6,-1.5)circle(0.1);
			\filldraw[fill=white](6.8,-1.5)circle(0.1);
			\filldraw[fill=white](7.6,-1.5)circle(0.1);
			\filldraw[fill=white](8.4,-1.5)circle(0.1);
			\filldraw[fill=white](9.2,-1.5)circle(0.1);
			
			\draw[very thick](0,0)--(-0.4,-1.5);
			\draw[very thick](0,0)--(0.4,-1.5);
			\draw[very thick](0.8,0)--(0.4,-1.5);
			\draw[very thick](0.8,0)--(1.2,-1.5);
			\draw[very thick](1.6,0)--(1.2,-1.5);
			\draw[very thick](1.6,0)--(2,-1.5);
			\draw[very thick](2.4,0)--(2,-1.5);
			\draw[very thick](2.4,0)--(2.8,-1.5);
			\draw[very thick](3.2,0)--(2.8,-1.5);
			\draw[very thick](3.2,0)--(3.6,-1.5);
			\draw[very thick](4,0)--(3.6,-1.5);
			\draw[very thick](4,0)--(4.4,-1.5);
			\draw[very thick](4.8,0)--(4.4,-1.5);
			\draw[very thick](4.8,0)--(5.2,-1.5);
			\draw[very thick](5.6,0)--(5.2,-1.5);
			\draw[very thick](5.6,0)--(6,-1.5);
			\draw[very thick](6.4,0)--(6,-1.5);
			\draw[very thick](6.4,0)--(6.8,-1.5);
			\draw[very thick](7.2,0)--(6.8,-1.5);
			\draw[very thick](7.2,0)--(7.6,-1.5);
			\draw[very thick](8,0)--(7.6,-1.5);
			\draw[very thick](8,0)--(8.4,-1.5);
			\draw[very thick](8.8,0)--(8.4,-1.5);
			\draw[very thick](8.8,0)--(9.2,-1.5);
			
			\draw[blue][dotted][very thick](2.8,0.5)--(3.2,-2);
			
			\coordinate [label=above:$29$] (x) at (0,0.1);
			\coordinate [label=above:$31$] (x) at (0.8,0.1);
			\coordinate [label=above:$33$] (x) at (1.6,0.1);
			\coordinate [label=above:$35$] (x) at (2.4,0.1);
			\coordinate [label=above:$37$] (x) at (3.2,0.1);
			\coordinate [label=above:$39$] (x) at (4,0.1);
			\coordinate [label=above:$41$] (x) at (4.8,0.1);
			\coordinate [label=above:$43$] (x) at (5.6,0.1);
			\coordinate [label=above:$45$] (x) at (6.4,0.1);
			\coordinate [label=above:$47$] (x) at (7.2,0.1);
			\coordinate [label=above:$49$] (x) at (8,0.1);
		    \coordinate [label=above:$51$] (x) at (8.8,0.1);
			\coordinate [label=above:$25$] (x) at (-0.4,-2.2);
			\coordinate [label=above:$23$] (x) at (0.4,-2.2);
			\coordinate [label=above:$21$] (x) at (1.2,-2.2);
			\coordinate [label=above:$19$] (x) at (2,-2.2);
			\coordinate [label=above:$17$] (x) at (2.8,-2.2);
			\coordinate [label=above:$15$] (x) at (3.6,-2.2);
			\coordinate [label=above:$13$] (x) at (4.4,-2.2);
			\coordinate [label=above:$11$] (x) at (5.2,-2.2);
			\coordinate [label=above:$9$] (x) at (6,-2.2);
			\coordinate [label=above:$7$] (x) at (6.8,-2.2);
			\coordinate [label=above:$5$] (x) at (7.6,-2.2);
			\coordinate [label=above:$3$] (x) at (8.4,-2.2);
			\coordinate [label=above:$1$] (x) at (9.2,-2.2);
		\end{tikzpicture}
	\end{center}
	\caption{The  representation of  $MD(\lambda)=\{51,47,43,39,35,31,25\}$.}\label{fig:lar2}
\end{figure}

For example, let $t=26$ and  $\lambda$ be a partition in $\mathcal{DS}(26)$ with the largest size. Then $51,47, 43, 39 \in MD(\lambda)$. Moreover, we have $\ell=\lfloor \frac{2\times26+1}{12}\rfloor+2=6$ and $\Delta=\frac{26}{2}-1-2\times 4=4$. It follows that $35,31,25\in MD(\lambda)$. Hence we obtain $MD(\lambda)=\{51,47,43,39,35,31,25\}$ as shown in Figure \ref{fig:lar2}. Then $|\lambda|=271$.

\noindent{\bf Case 3.} When $\Delta=\frac{t}{2}-1-2\lfloor \frac{2t+1}{12}\rfloor \equiv  2 \, (\text{mod}\ 3)$. Similarly, combining  Lemmas \ref{lle2}, \ref{llea}, \ref{lle3}, and \ref{lle4}, we derive that  $t+5, t-1 \in MD(\lambda)$ when $\ell \ne k$.
% recall that $t+5, t-1 \in MD(\lambda)$ when $\ell \ne k$. 
Furthermore,  we derive that $ 2t-1-4\lfloor \frac{2t+1}{12}\rfloor, 2t-7-4\lfloor \frac{2t+1}{12}\rfloor, \dots, t+5,t-1,\dots, 5+4\lfloor \frac{2t+1}{12}\rfloor \in MD(\lambda)$. Recall that $2t+3-4\lfloor \frac{2t+1}{12}\rfloor \in MD(\lambda)$. Then we deduce that $\ell=\lfloor \frac{2t+1}{12}\rfloor+1$, namely, $d(\lambda)=(d_1,\dots,d_{\lfloor \frac{2t+1}{12}\rfloor},d_{\lfloor \frac{2t+1}{12}\rfloor+1},\dots, d_{k-1})=(4,\dots,4,6,\dots,6)$. Thus,  we have
$$
MD(\lambda)=\{2t-1,2t-5,\dots,2t-1-4\lfloor \frac{2t+1}{12}\rfloor, 2t-7-4\lfloor \frac{2t+1}{12}\rfloor, \dots, 5+4\lfloor \frac{2t+1}{12}\rfloor \}.
$$
According to Lemma \ref{sclar}, the size of $\lambda$ is 
$$
|\lambda|
=\frac{1}{3}(t-1)\left(t-3-4\lfloor \frac{2t+1}{12} \rfloor \right)+\left(\lfloor \frac{2t+1}{12} \rfloor+1\right) \left(2t-1-2\lfloor \frac{2t+1}{12} \rfloor \right)
$$
as desire.

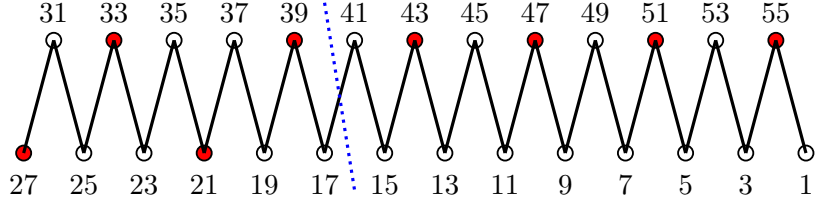
\begin{figure}[H]
	\begin{center}
		\begin{tikzpicture}[font =\small , scale = 1, line width = 0.7pt]
			\filldraw[fill=white](0,0)circle(0.1);
			\filldraw[fill=red](0.8,0)circle(0.1);
			\filldraw[fill=white](1.6,0)circle(0.1);
			\filldraw[fill=white](2.4,0)circle(0.1);
			\filldraw[fill=red](3.2,0)circle(0.1);
			\filldraw[fill=white](4,0)circle(0.1);
			\filldraw[fill=red](4.8,0)circle(0.1);
			\filldraw[fill=white](5.6,0)circle(0.1);
			\filldraw[fill=red](6.4,0)circle(0.1);
			\filldraw[fill=white](7.2,0)circle(0.1);
			\filldraw[fill=red](8,0)circle(0.1);
			\filldraw[fill=white](8.8,0)circle(0.1);
			\filldraw[fill=red](9.6,0)circle(0.1);
			
			\filldraw[fill=red](-0.4,-1.5)circle(0.1);
			\filldraw[fill=white](0.4,-1.5)circle(0.1);
			\filldraw[fill=white](1.2,-1.5)circle(0.1);
			\filldraw[fill=red](2,-1.5)circle(0.1);
			\filldraw[fill=white](2.8,-1.5)circle(0.1);
			\filldraw[fill=white](3.6,-1.5)circle(0.1);
			\filldraw[fill=white](4.4,-1.5)circle(0.1);
			\filldraw[fill=white](5.2,-1.5)circle(0.1);
			\filldraw[fill=white](6,-1.5)circle(0.1);
			\filldraw[fill=white](6.8,-1.5)circle(0.1);
			\filldraw[fill=white](7.6,-1.5)circle(0.1);
			\filldraw[fill=white](8.4,-1.5)circle(0.1);
			\filldraw[fill=white](9.2,-1.5)circle(0.1);
			\filldraw[fill=white](10,-1.5)circle(0.1);
			
			\draw[very thick](0,0)--(-0.4,-1.5);
			\draw[very thick](0,0)--(0.4,-1.5);
			\draw[very thick](0.8,0)--(0.4,-1.5);
			\draw[very thick](0.8,0)--(1.2,-1.5);
			\draw[very thick](1.6,0)--(1.2,-1.5);
			\draw[very thick](1.6,0)--(2,-1.5);
			\draw[very thick](2.4,0)--(2,-1.5);
			\draw[very thick](2.4,0)--(2.8,-1.5);
			\draw[very thick](3.2,0)--(2.8,-1.5);
			\draw[very thick](3.2,0)--(3.6,-1.5);
			\draw[very thick](4,0)--(3.6,-1.5);
			\draw[very thick](4,0)--(4.4,-1.5);
			\draw[very thick](4.8,0)--(4.4,-1.5);
			\draw[very thick](4.8,0)--(5.2,-1.5);
			\draw[very thick](5.6,0)--(5.2,-1.5);
			\draw[very thick](5.6,0)--(6,-1.5);
			\draw[very thick](6.4,0)--(6,-1.5);
			\draw[very thick](6.4,0)--(6.8,-1.5);
			\draw[very thick](7.2,0)--(6.8,-1.5);
			\draw[very thick](7.2,0)--(7.6,-1.5);
			\draw[very thick](8,0)--(7.6,-1.5);
			\draw[very thick](8,0)--(8.4,-1.5);
			\draw[very thick](8.8,0)--(8.4,-1.5);
			\draw[very thick](8.8,0)--(9.2,-1.5);
			\draw[very thick](9.6,0)--(9.2,-1.5);
			\draw[very thick](9.6,0)--(10,-1.5);
			
		\draw[blue][dotted][very thick](3.6,0.5)--(4,-2);
			
			\coordinate [label=above:$31$] (x) at (0,0.1);
			\coordinate [label=above:$33$] (x) at (0.8,0.1);
			\coordinate [label=above:$35$] (x) at (1.6,0.1);
			\coordinate [label=above:$37$] (x) at (2.4,0.1);
			\coordinate [label=above:$39$] (x) at (3.2,0.1);
			\coordinate [label=above:$41$] (x) at (4,0.1);
			\coordinate [label=above:$43$] (x) at (4.8,0.1);
			\coordinate [label=above:$45$] (x) at (5.6,0.1);
			\coordinate [label=above:$47$] (x) at (6.4,0.1);
			\coordinate [label=above:$49$] (x) at (7.2,0.1);
			\coordinate [label=above:$51$] (x) at (8,0.1);
			\coordinate [label=above:$53$] (x) at (8.8,0.1);
			\coordinate [label=above:$55$] (x) at (9.6,0.1);
			\coordinate [label=above:$27$] (x) at (-0.4,-2.2);
			\coordinate [label=above:$25$] (x) at (0.4,-2.2);
			\coordinate [label=above:$23$] (x) at (1.2,-2.2);
			\coordinate [label=above:$21$] (x) at (2,-2.2);
			\coordinate [label=above:$19$] (x) at (2.8,-2.2);
			\coordinate [label=above:$17$] (x) at (3.6,-2.2);
			\coordinate [label=above:$15$] (x) at (4.4,-2.2);
			\coordinate [label=above:$13$] (x) at (5.2,-2.2);
			\coordinate [label=above:$11$] (x) at (6,-2.2);
			\coordinate [label=above:$9$] (x) at (6.8,-2.2);
			\coordinate [label=above:$7$] (x) at (7.6,-2.2);
			\coordinate [label=above:$5$] (x) at (8.4,-2.2);
			\coordinate [label=above:$3$] (x) at (9.2,-2.2);
			\coordinate [label=above:$1$] (x) at (10,-2.2);
		\end{tikzpicture}
	\end{center}
	\caption{The  representation of  $MD(\lambda)=\{55,51,47,43,39,33,27,21\}$.}\label{fig:lar3}
\end{figure}

For example, let $t=28$ and  $\lambda$ be a partition in $\mathcal{DS}(28)$ with the largest size. Then $55, 51, 47, 43 \in MD(\lambda)$, Moreover, we have $\ell=\lfloor \frac{2\times28+1}{12}\rfloor+1=5$ and $\Delta=\frac{28}{2}-1-2\times 4=5$. It follows that $39,33,27,21 \in MD(\lambda)$. Hence we obtain  $MD(\lambda)=\{55,51,47,43,39,33,27,21\}$ as shown in Figure \ref{fig:lar3}. Then $|\lambda|=316$.

When $t$ is odd and $t \geq 6$. 
% Combining  Lemma \ref{lle2},  \ref{llea},  \ref{lle3}, and  \ref{lle4}, we derive that  $t-4, t+2 \in MD(\lambda)$ when $\ell \ne k$.  On the other hand, for a given odd positive integer $t$, i
It is not difficult to check that $|\mathrm{Q}_{t}^{l}|=\frac{t-1}{2}$. Recall that $\{2t-1, 2t-5, \dots, 2t+3-4\lfloor \frac{2t+1}{12}\rfloor\} \subseteq MD(\lambda)$. Then we have $2t-3, 2t-7, \dots, 2t+1-4\lfloor \frac{2t+1}{12}\rfloor \notin MD(\lambda)$ by the definition of partitions in $\mathcal{DS}(t)$. Thus we obtain that
$$|\mathrm{Q}_{t}^{l}\backslash \{2t-1,2t-3,\dots, 2t+1-4\lfloor \frac{2t+1}{12}\rfloor\}|=\frac{t-1}{2}-2\lfloor \frac{2t+1}{12}\rfloor.$$ 

Similarly, we have three cases:  {\upshape (\rmnum{1})}  $\frac{t-1}{2}-2\lfloor \frac{2t+1}{12}\rfloor \equiv 0 \ (\text{mod}\ 3)$;  {\upshape (\rmnum{2})}  $\frac{t-1}{2}-2\lfloor \frac{2t+1}{12}\rfloor \equiv 1 \ (\text{mod}\ 3)$; and {\upshape (\rmnum{3})}  $\frac{t-1}{2}-2\lfloor \frac{2t+1}{12}\rfloor \equiv 2 \ (\text{mod} \ 3)$. 
All these cases can be verified by similar arguments to the case $t \geq 6$ be odd and will be omitted, completing the proof of Theorem \ref{largest}.\qed

\section*{Acknowledgments}

The authors are grateful to the anonymous referees for valuable comments and suggestions which helped improve the paper. This work was supported by the National Science Foundation of China [Grant No. 12201155].

\end{document}